\documentclass[twocolumn]{autart}    

\usepackage{graphicx}      

\usepackage{amssymb}  
\usepackage[shortlabels]{enumitem}  
\usepackage{array}
\usepackage{nicematrix}
\usepackage{color}
\usepackage{pifont}
\usepackage{tabularx} 
\usepackage{balance}
\usepackage{algorithm}

\usepackage{algorithmic}
\usepackage{multirow}
\usepackage{makecell}
\usepackage{bbm}

\begin{document}

\begin{frontmatter}

\title{A differential game approach to intrinsic encirclement control\thanksref{footnoteinfo}} 

\thanks[footnoteinfo]{This paper is supported by KTH Digital Futures.}

\author[First]{Panpan Zhou} \ead{panpanz@kth.se},
\author[Second]{Yueyue Xu}\ead{merryspread99@sjtu.edu.cn}, 
\author[Third]{Yibei Li}\ead{yibeili@amss.ac.cn}, 
\author[Fourth]{Bo Wahlberg} \ead{bo@kth.se},
\author[First]{Xiaoming Hu*} \ead{hu@kth.se}


\address[First]{Department of Mathematics, KTH Royal Institute of Technology, Stockholm, Sweden}
\address[Second]{Department of Automation, Shanghai Jiao Tong University, Shanghai 201100, China}
\address[Third]{Key Laboratory of Systems and Control, Academy of Mathematics and Systems Science,
        Chinese Academy of Sciences, Beijing 100190, P.~R.~China}
\address[Fourth]{Division of Decision and Control Systems, KTH Royal Institute of Technology, Stockholm, Sweden}



\begin{abstract}                          
This paper investigates the encirclement control problem involving two groups using a non-cooperative differential game approach. The active group seeks to chase and encircle the passive group, while the passive group responds by fleeing cooperatively and simultaneously encircling the active group. Instead of prescribing an expected radius or a predefined path for encirclement, we focus on the whole formation manifold of the desired relative configuration, two concentric circles, by allowing permutation, rotation, and translation of players. The desired relative configurations arise as the steady state resulting from Nash equilibrium strategies and are achieved in an intrinsic way by designing the interaction graphs and weight function of each edge.
Furthermore, the asymptotic convergence to the desired manifold is guaranteed. Finally, numerical simulations demonstrate encirclement and counter-encirclement scenarios, verifying the effectiveness of our strategies.
\end{abstract}

\end{frontmatter}

\section{Introduction}\label{sec1}

Encirclement control is one of the most significant motion control problems in the field of multi-agent systems, presenting both theoretical challenges and practical potential across various applications, such as environmental surveys \cite{slingsby2023review} and source seeking \cite{han2014multiple}. Closely related fields include containment control \cite{zhou2020formation,cheng2023data,huang2022distributed,zhou2023distributed} and fencing control \cite{kou2021cooperative1,kou2021cooperative2}, all of which aim to enclose targets while maintaining specific spatial relationships among agents. However, encirclement control specifically emphasizes circular motion or configuration around targets.

Generally, most existing methods define an expected distance or position around the target as a reference and then regulate the enclosing error. For instance, in encirclement control, references \cite{liu2023multiple,yue2021elliptical} defined an expected constant or time-varying radius, and reference \cite{zhang2024optimal} specified a desired path around the target. Similarly, in target fencing \cite{hu2021distributed}, a predefined distance between agents and the target was used as a control references. However, in many application scenarios, obtaining the absolute enclosing error for each agent is challenging. 
Instead, we address the encirclement control while maintaining formation patterns by focusing on the relative configuration of the entire multi-agent system, allowing for possible permutation, orientation and translation of the formation. Therefore, in this paper, we solve the encirclement control problem in an intrinsic way, meaning that the desired formation emerges solely from inter-agent interactions and the geometric properties of the network, without predefined desired formation or fixed positions assigned to each agent in the controller. The concept of intrinsic control originates from \cite{song2017intrinsic}, where reduced attitude formation is studied in $\mathcal{S}^2$. Artificial potential functions have also been employed in \cite{fabiani2018distributed,ogren2004cooperative} to achieve families of equilibrium solutions exhibiting both rotational and translational symmetry.

Until now, most studies on encirclement and fencing control focus on a single target \cite{chen2019cooperative,hu2021distributed,zhang2024elliptical,shevchenko2017guaranteed,zhang2024optimal}. Although multiple targets are considered in some works, such as \cite{hafez2013encirclement,jiang2022constrained,hafez2013encirclement}, targets are either stationary \cite{hafez2013encirclement} or follow fixed trajectories \cite{hafez2013encirclement,jiang2022constrained}, with an emphasis on estimating their positions. The containment control is typically studied in the leader-follower framework, where the objective is to drive the followers toward the convex hull formed by the leaders \cite{gonzalez2024saturated,cheng2023data}, even when the leaders follow dynamic trajectories \cite{zhou2020formation}. 
It should be noted that a common feature of these results is that the target(s) or leaders can not respond to the actions of other agents. However, in more practical scenarios, targets may evade others and even attempt counter-encirclement. In this paper, we study both encirclement control and counter-encirclement control with considering interactions between two groups.

In recent years, game-theoretic approaches have provided fresh insights into control problems, especially in formation control \cite{li2022differential,fang2022formation}, where agents either cooperate to complete a shared task or optimize similar objective functions. However, results are still lacking for more practical scenarios involving groups with individual and conflicting objectives, leading to a non-cooperative setting. Within the game-theoretic framework, each agent seeks to find its Nash Equilibrium strategy to minimize its own cost while accounting for the strategies of other agents. The game is theoretically challenging even with quadratic functions, let alone dealing with non-quadratic objective functions. In this paper, we focus on the game with non-quadratic functions. 
Different from existing formation control approaches based on differential games theory, which typically formulate quadratic games and regulate the formation error, the game with intrinsic property is outside the consensus framework. In addition to the intrinsic control in \cite{song2017intrinsic}, the work \cite{li2022differential} studies the intrinsic formation control for more general patterns in $\mathcal{R}^3$, where Nash equilibrium strategies ensure convergence to the desired formation manifold.

The contributions of this paper are as follows. 
First, compared to previous works that directly design control laws \cite{santilli2022secure,gonzalez2024saturated,liu2023multiple}, we solve the encirclement control problem from a non-cooperative differential game perspective. This brings theoretical challenges in obtaining Nash Equilibrium strategies, as the problem involves non-quadratic objective functions within a non-consensus framework.
Second, we investigate encirclement control between two groups of agents with opposing objectives, inspired by the natural predator-prey dynamics of encirclement and counter-encirclement. To the best of our knowledge, studies on mutual encirclement are rare. Our work differs from most existing literature, where the targets do not respond to the actions of other agents \cite{chen2019cooperative,shevchenko2017guaranteed,gonzalez2024saturated}.
Third, instead of pre-assigning desired positions for agents, we achieve the desired formation manifold solely by designing the graph topology and the weight functions of its edges. This approach presents additional challenges, as it requires accommodating greater flexibility in defining the desired configuration.
Fourth, the desired relative pattern is shown to be asymptotic stable. The convergence follows from the invariant covariant tensor field and symmetric properties of edge weights and the desired formation pattern.

The rest of this paper is organized as follows. In Section~\ref{sec2}, the encirclement control problem is formulated as a non-cooperative differential game.  Section~\ref{sec3} presents properties preserved under rotation and translation group actions, along with the designed topologies. Based on these properties, the infinite time-horizon game is analyzed in Section~\ref{sec4}. In Section~\ref{sec5} numerical simulations are provided to validate the effectiveness of our strategies. Some concluding remarks are given in Section~\ref{sec6}.

\textit{Notations:} $\|\cdot\|$ denotes the $l_2$ norm of a vector. $I_n\in \mathbb{R}^{n\times n}$ is an identity matrix. $\mathbbm{1}_m:=[1~1~\cdots~1]^{\rm T}\in \mathbb{R}^m$ is a vector whose elements are zero. ${\rm col}\{x,y\}:=[x~y]^{\rm T}$. $\otimes$ denotes the Kronecker product. ${\rm card}(\cdot)$ represents the cardinality of a set. $diag(a_1,a_2,\cdots,a_m)$ denotes a matrix whose diagonal elements are $a_1$, $a_2$, $\cdots$, and all non-diagonal elements are zero. $\lambda_k(\cdot)$ represents the $k$-th eigenvalue of a matrix.

\section{Problem formulation}\label{sec2}
In this work, we consider two groups, referred to as the active group and the passive group, which play a non-cooperative differential game with each group consisting of $m$ and $n$ players, respectively. The active group chases and encircles the passive group, while the passive group responds by cooperatively fleeing and also encircling the active group. The dynamics of each player is modeled as a single integrator, i.e.,
\begin{align}
    \dot x_i^a &= u_i^a, \quad i=1,2,\cdots,m, \label{group_a}\\
    \dot x_j^p &= u_j^p, \quad j=1,2,\cdots,n, \label{group_p}
\end{align}
where $x_i^a,x_j^p\in \mathbb{R}^2$ are the positions of the $i$th active player and the $j$th passive player, respectively. $u_i^a,u_j^p$ denote their respective control inputs. 
Let $x={\rm col}\{x_1^a,\cdots,x_m^a,x_1^p,\cdots,x_n^p\}$. The full-state dynamics of the multi-agent system are given by
\begin{equation}  \label{dyn}
    \dot x=\sum_{l=1}^{m+n}B_lu_l 
\end{equation}
where $B_l\in \mathbb{R}^{2(m+n)\times 2}$ is a column-wise block matrix consisting of $m+n$ blocks of size $2\times 2$, with only the $l$th block being an identity matrix and all other blocks being zero matrices.
For $l=1,2,\cdots,m$, $u_l=u_l^a$, while for $l=m+1,m+2,\cdots,m+n$, $u_l=u_{l-m}^p$.

We assume that each player exchanges information with its connected neighbors through a communication graph $\mathcal{G}=\mathcal{G}^a \cup \mathcal{G}^p$. Specifically, $\mathcal{G}^a=(\mathcal{V}^a,\mathcal{E}^a)$ where $\mathcal{V}^a=\{1,2,\cdots,m\}$ denotes the vertex set of the active group, and $\mathcal{E}^a \subset\mathcal{V}^a\times \mathcal{V}^a$ is the edge set. $(k,i)\in \mathcal{E}^a$ if and only if player $i$ has access to the information of player $k$, and we say player $k$ is a neighbor of player $i$. Similarly, we define $\mathcal{G}^p=(\mathcal{V}^p,\mathcal{E}^p)$ for the passive group. The set of neighbors of player $i$ is denoted by $\mathcal{N}_i=\{k\in \mathcal{V}:(k,i)\in \mathcal{E}\}$. For $i\in \mathcal{V}^a$, we define $\mathcal{N}_i^a$ as the subset of its neighbors belonging to the active group, and let $d_i^a={\rm card}(\mathcal{N}_i^a)$ denote the number of such neighbors. Similarly, for $j\in \mathcal{V}^p$, we define $\mathcal{N}_j^p$ and $d_j^p$ to represent, respectively, the set of its neighbors and the number of its neighbors in the passive group. 
In this work, we consider \textit{undirected graphs}, namely, $(k,i)\in \mathcal{E}$ means $(i,k)\in \mathcal{E}$. We also assume that every member of the active group has access to the positions of all members in the passive group, and vice versa.

The encirclement control problem of the multi-agent system is modeled as a non-cooperative game, where each player is associated with an individual cost function of the form 
\begin{align}
    J_i^a(x(0),u_i^a,u_{-i})=\frac{1}{2} \int_0^{\infty} (q_i^a(x)+\|u_i^a\|^2)dt,  \label{obj_a} \\
    J_j^p(x(0),u_j^p,u_{-j})=\frac{1}{2} \int_0^{\infty} (q_j^p(x)+\|u_j^p\|^2)dt,   \label{obj_p}
\end{align}
for $i\in \mathcal{V}^a$ and $j\in \mathcal{V}^p$, 
where $u_{-i}$ is the strategy profile of all players except player $i$. Combining with the full-state dynamics (\ref{dyn}), the non-cooperative game played by the $i$th active player is formulated as
\begin{equation} \label{game_a}
    \begin{aligned}
        &\min_{u_i^a}~~J_i^a(x(0),u_i^a,u_{-i}) \\
        &s.t. ~~~~\dot x=\sum_{l=1}^{m+n}B_lu_l .
    \end{aligned}
\end{equation}
Similarly, the non-cooperative game played by the $j$th passive player is formulated as
\begin{equation} \label{game_p}
    \begin{aligned}
        &\min_{u_j^p}~~J_j^p(x(0),u_j^p,u_{-j})\\
        &s.t. ~~~~\dot x=\sum_{l=1}^{m+n}B_lu_l .
    \end{aligned}
\end{equation}

Before defining the problem, we provide definitions of some Nash equilibrium-related concepts.

\newtheorem{definition}{Definition}
\begin{definition}
    (Nash equilibrium and Nash equilibrium strategies) Consider the non-cooperative games (\ref{game_a}) and (\ref{game_p}). A strategy profile $u^*=(u_1^{a*},\cdots,u_m^{a*},u_1^{p*},\cdots,u_n^{p*})$ is a Nash equilibrium strategy if the following conditions hold:
    \begin{equation*}
    \begin{aligned}
        J_i^a(x(0),u_i^{a*},u_{-i}^*) \le J_i^a(x(0),u_i^a,u_{-i}^*), ~i\in \mathcal{V}^a,\\
        J_j^p(x(0),u_j^{p*},u_{-j}^*) \le J_i^p(x(0),u_j^p,u_{-j}^*), ~j\in \mathcal{V}^p,
    \end{aligned}
    \end{equation*}
    where $u_{-i}^*=(u_1^{a*},\cdots,u_{i-1}^{a*},u_{i+1}^{a*},\cdots, u_m^{a*},u_1^{p*},\cdots,u_n^{p*})$ and $u_{-j}^*=(u_1^{a*},\cdots,u_m^{a*},u_1^{p*},\cdots,u_{j-1}^{p*},u_{j+1}^{p*},\cdots,u_n^{p*})$ represent the strategy profiles of all players except player $i$ and $j$, respectively. The steady state resulting from Nash equilibrium strategies is referred to as the Nash equilibrium itself. Besides, the trajectory $x$ under $u^*$ is called the Nash equilibrium trajectory.
\end{definition}

The solution $u^*$ is referred to as a Nash equilibrium in most literature on game theory \cite{isaacs1999differential,bacsar1998dynamic,de2017finite} as they focused on solutions where players cannot improve their objectives by deviating from it, without concern for the final outcome of the game. In contrast, this paper values not only the solution but also prioritizes the steady state. Accordingly, we distinguish between Nash equilibrium strategies (the solution $u_i^*$) and the Nash equilibrium (the resulting steady state). This distinction aligns with our objective of guiding players toward a desired configuration under criteria (\ref{game_a}) and (\ref{game_p}).

In the following, we introduce some similar preliminary results on permutation, rotation, and translation as in \cite{li2022differential}. We deem this introduction necessary since in this paper we do not consider the permutations between the active and passive players, whereas \cite{li2022differential} allows permutations among all players; and we focus on rotation in 2D rather than 3D. 

A coordination-based description is introduced to characterize the manifold of desired vertex coordinates. An undirected graph $\mathcal{P}$ is used to denote the target relative configuration with its edges defining the skeletal structure. Let $\mathcal{S}_{\mathcal{P}}^a$ and $\mathcal{S}_{\mathcal{P}}^p$ denote the permutation groups for the active group and the passive group, respectively, corresponding to the sets of bijections $\sigma^a:\mathcal{V}^a\to \mathcal{V}^a$ and $\sigma^p:\mathcal{V}^p\to \mathcal{V}^p$ under function composition. Each permutation $\sigma^a\in \mathcal{S}_{\mathcal{P}}^a$ and $\sigma^p\in \mathcal{S}_{\mathcal{P}}^p$ can be represented by permutation matrices $P_{\sigma^a}=[e_{\sigma^a(1)},\cdots,e_{\sigma^a(m)}]$ and $P_{\sigma^p}=[e_{\sigma^p(1)},\cdots,e_{\sigma^p(n)}]$, where $e_{\sigma(k)}$ is a column vector with 1 at element $\sigma(k)$ and 0 elsewhere. In this paper, the active players cannot permute with the passive players, as they have different objectives. Let $x^*=[x_a^{* \rm T}~x_p^{* \rm T}]^{\rm T}={\rm col}\{x_1^{a*},\cdots,x_m^{a*},x_1^{p*},\cdots,x_n^{p*}\}$ be the vertex coordinates of the desired pattern centered at the origin. Then, through vertex permutation and rotation of the body, the whole vertices set of $\mathcal{P}$ is given by
\begin{equation}
\begin{aligned}
    \mathcal{M}_{\mathcal{P}}\!&=\! \{x\!\in\!\mathbb{R}^{2(m+n)}\!:\!x\!=\!(I_{m+n} \!\otimes\! R)(diag(P_{\sigma^a}\!,P_{\sigma^p})\!\otimes\! I_2)x^*, \\
    &\forall R\in SO(2), \sigma^a \in \mathcal{S}_{\mathcal{P}}^a, \sigma^p \in \mathcal{S}_{\mathcal{P}}^p\}.
\end{aligned}
\end{equation}

Let $o^a$ represent the order of the group $\mathcal{S}_{\mathcal{P}}^a$ with elements $\{\sigma_1^a,\cdots,\sigma_{o^a}^a\}$. Similarly, let $o^p$ denote the order of the group $\mathcal{S}_{\mathcal{P}}^p$ with elements $\{\sigma_1^p,\cdots,\sigma_{o^p}^p\}$. Then, the manifold $\mathcal{M}_{\mathcal{P}}$ can be expressed as the union of a finite number of disjoint manifolds, i.e., $\mathcal{M}_{\mathcal{P}}= \bigcup  \mathcal{M}_{kl}$ where
\begin{equation}
\begin{aligned} \label{manif}
    \mathcal{M}_{kl}\!=\!&\{x \!\in\! \mathbb{R}^{2(m+n)}\!:\! x\!=\!(I_{m+n} \!\otimes\! R)(diag(P_{\sigma_k^a}\!,P_{\sigma_l^p})\!\otimes\! I_2)x^*, \\
     &\forall R\in SO(2), k=1,\cdots,o^a, l=1,\cdots,o^p.
\end{aligned}
\end{equation}

\newtheorem{proposition}{Proposition}
\begin{proposition}
    Each manifold $\mathcal{M}_{kl}$ is closed. For any $k_1,k_2=1,\cdots,o^a$, $l_1,l_2=1,\cdots,o^p$, $\mathcal{M}_{k_1,l_1}$, $\mathcal{M}_{k_1,l_2}$, $\mathcal{M}_{k_2,l_1}$ and $\mathcal{M}_{k_2,l_2}$ are either equal or disjoint.
\end{proposition}
The proof follows a similar approach to that of Proposition 2.1 in \cite{li2022differential}.

Each manifold $\mathcal{M}_{kl}$ in (\ref{manif}) is obtained by rotating the configuration that corresponds to a specific vertex ordering. Since players are engaged in chasing and fleeing in 2D, more degrees of freedom are introduced to characterize the vertex manifold of the desired configuration, i.e., two concentric circles with a translatable center.
We define the \textit{translation manifold} by successively applying the rotation and translation groups to $x^*$ as
\begin{equation*}
\begin{aligned}
    \mathcal{M}_{\rm T}(x^*)=\{x \!\in\! \mathbb{R}^{2(m\!+\!n)}:x&\!=\!(I_{m+n}\otimes R)x^* \!+\! \mathbbm{1}_{m\!+\!n}\otimes s, \\
    & \forall R\in SO(2), \forall s\in \mathbb{R}^2\}.
\end{aligned}
\end{equation*}

Based on the above information, we give the definition of the problem.

\newtheorem{problem}{Problem}
\begin{problem}
    Consider a multi-agent system playing a non-cooperative game as defined in (\ref{game_a})--(\ref{game_p}). Design $q_i^a(x)$ and $\mathcal{N}_i^a$ for the active group, and $q_j^p(x)$ and $\mathcal{N}_j^p$ for the passive group, such that the resulting Nash equilibrium trajectory locally converges to the desired pattern $\mathcal{P}$ where one group encircles the other. It means that for any $i\in \mathcal{V}^a$ and $j\in \mathcal{V}^p$, the Nash equilibrium strategies $(u_1^{a*},\cdots,u_m^{a*}, u_1^{p*},\cdots,u_n^{p*})$ satisfy
    \begin{equation*}
    \begin{aligned}
        J_i^a(x(0),u_i^{a*},u_{-i}^*) \le J_i^a(x(0),u_i^a,u_{-i}^*), ~i\in \mathcal{V}^a,\\
        J_j^p(x(0),u_j^{p*},u_{-j}^*) \le J_j^p(x(0),u_j^p,u_{-j}^*), ~j\in \mathcal{V}^p,
    \end{aligned}
    \end{equation*}
    for all $u_i^a\ne u_i^{a*}$ and $u_j^p\ne u_j^{p*}$, and $\mathcal{M}_{\rm T}(x^*)$ is asymptotically stable under Nash equilibrium strategies.
\end{problem}

\section{Symmetry and invariance on manifolds} \label{sec3}
In this paper, we consider a formation where two concentric circles are formed in 2D, with each group forming one circle. The players are evenly distributed on their respective circle. It can be seen that when the number of players on the circle is even, the configuration is highly symmetric.

\subsection{Properties of rotation group and translation group actions}\label{sec3.1}
In this subsection, we present some preliminary results that highlight the properties preserved under rotational symmetry and translation.

The tangent space at point $x\in \mathcal{M}_{\rm T}(x^*)$ is given by
\begin{equation} \label{mm3}
    \mathcal{T}_{\rm T}(x^*)\!=\!\underbrace{\{(I_{m\!+\!n}\otimes L)x^*:L \!\in\! \mathfrak{so}(2)\}}_{\mathcal{T}_R(x^*)} \!+\! \underbrace{\{\mathbbm{1}_{m\!+\!n}\otimes v:v \!\in \! \mathbb{R}^2\} }_{\mathcal{T}_v}.
\end{equation}

\begin{proposition}
For any non-collinear coordinate of vertices $x^*$, ${\rm dim}(\mathcal{T}_{\rm T}(x^*))=3$.
\end{proposition} 
\textbf{Proof.} It is evident that ${\rm dim}(\mathcal{T}_v)=2$ and ${\rm dim}(\mathcal{T}_R(x^*))=1$ for any non-collinear $x^*$. Then we demonstrate that $\mathcal{T}_v$ and $\mathcal{T}_R(x^*)$ are linearly independent by contradiction. Assume  there exists $\bar x={\rm col}\{\bar x_1,\cdots,\bar x_{m+n}\}\in \mathcal{T}_R(x^*) \cap \mathcal{T}_v$ with $\bar x\ne 0$. Then there must exist a nonzero vector $u_L\in\mathbb{R}^2$ such that $u_L\times \bar x_i=u_L\times \bar x_j$ for any $i,j\in \mathcal{V}$, which contradicts the assumption that $\{\bar x\}_1^{m+n}$ are non-collinear.

\begin{definition}
    ($\mathcal{R}$-invariant). A real-valued function $V(x)$ is called $\mathcal{R}$-invariant if it is invariant under rotation groups, i.e., $V(\bar Rx)=V(x)$, for any $\bar R\in SO(2)$. 
\end{definition}

It can be proven that if $V(x)$ is $\mathcal{R}$-invariant, its gradient is also invariant in the sense that $\frac{\partial V(x)}{\partial x}|_{x=\bar Rx}=\frac{\partial V(x)}{\partial x}\bar R^{\rm T}$. Furthermore, the Hessian of $V(x)$ satisfies $\frac{\partial^2 V(x)}{\partial x^2}|_{x=\bar Rx}=\bar R \frac{\partial^2 V(x)}{\partial x^2} \bar R^{\rm T}$.

\begin{definition}
    (Translation invariant). A real-valued function $V(x):\mathbb{R}^{2(m+n)}\to \mathbb{R}$ is called translation invariant if $V(x+\mathbbm{1}_{m+n}\otimes s)=V(x)$ for any $s\in \mathbb{R}^2$.
\end{definition}

The intuition behind this is that for any desired pattern, constructing its rotation and translation manifold as an invariant flow ensures that the properties of the entire manifold can be preserved by evaluating only an arbitrary point on it. This significantly simplifies the design and analysis of the corresponding games.
Then, we study some properties of the rotation and translation manifold.

\newtheorem{lemma}{Lemma}
\begin{lemma} \label{lemma1}
    \cite{li2022differential} Consider a vector field $f(x):\mathbb{R}^{2(m+n)}\to \mathbb{R}^{2(m+n)}$ that is $\mathcal{R}$-invariant. Then, for any equilibrium $x^*$, i.e., $f(x^*)=0$, it holds that
    \begin{equation*}
        \mathcal{T}_R(x^*)\in {\rm Ker}(\frac{\partial f(x^*)}{\partial x}).
    \end{equation*}
\end{lemma}

\newtheorem{theorem}{Theorem}
\begin{theorem} \label{them1}
    Consider $V(x):\mathbb{R}^{2(m+n)}\to \mathbb{R}$ that is $\mathcal{R}$-invariant and translation invariant. If there exists $x^*=[x_a^{* \rm T}~x_p^{* \rm T}]^{\rm T}\in \mathbb{R}^{2(m+n)}$ such that $\frac{\partial V(x^*)}{\partial x}=0$ and $\frac{\partial^2 V(x^*)}{\partial x^2}\ge 0$ with the number of zero eigenvalues equal to ${\rm dim}(\mathcal{M}_{\rm T}(x^*))$, then 
    $${\rm Ker}(\frac{\partial^2 V(x^*)}{\partial x^2})=\mathcal{T}_{\rm T}(x^*),$$
    and $V(x)$ obtains a strict local minimum at $\mathcal{M}_{\rm T}(x^*)$. Moreover, write $\frac{\partial^2 V(x^*)}{\partial x^2}$ as the form 
    \begin{equation} \label{mm4}
        \frac{\partial^{2} V(x^{*})}{\partial x^{2}}=\left[\begin{array}{c|c}
    A_1(x_{a}^{*},x_{p}^{*}) & A_{2}(x_{a}^{*},x_{p}^{*}) \\
    \hline
    A_{2}^{\rm T} (x_{a}^{*},x_{p}^{*}) & A_{3}(x_{a}^{*},x_{p}^{*})
    \end{array}\right]
    \end{equation}
    where $A_1(x_{a}^{*},x_{p}^{*})\in \mathbb{R}^{2m\times 2m}$, $A_{2}(x_{a}^{*},x_{p}^{*})\in \mathbb{R}^{2m\times 2n}$, $A_{3}(x_{a}^{*},x_{p}^{*})\in \mathbb{R}^{2n\times 2n}$. 
    If 
    \begin{equation}   \label{mm1}
        A_2(x_{a}^{*},x_{p}^{*})(I_n\otimes L)x_p^*=0, \\
    \end{equation}
    where $L$ is defined in (\ref{mm3}), 
    then the number of zero eigenvalues of $\frac{\partial^2 V(x^*)}{\partial x^2}$ equals to  ${\rm dim}(\mathcal{T}_{\rm T}(x^*))$ if and only if the number of zero eigenvalues of $A_1(x_a^*,x_p^*)$ equals to $\mathcal{M}_R(x_a^*)$, where $\mathcal{M}_R(x_a^*)=\{x\in \mathbb{R}^{2m}:x=(I_m\otimes R)x_a^*,\forall R\in SO(2)\}$; or if
    \begin{equation} \label{mm6}
        A_2^{\rm T}(x_{a}^{*},x_{p}^{*})(I_n\otimes L)x_a^*=0, \\
    \end{equation}
    then the number of zero eigenvalues of $\frac{\partial^2 V(x^*)}{\partial x^2}$ equals to ${\rm dim}(\mathcal{M}_{\rm T}(x^*))$ if and only if the number of zero eigenvalues of $A_3(x_a^*,x_p^*)$ equals to ${\rm dim}(\mathcal{M}_R(x_p^*))$, where $\mathcal{M}_R(x_p^*)=\{x\in \mathbb{R}^{2n}:x=(I_n\otimes R)x_p^*,\forall R\in SO(2)\}$.
\end{theorem}

\textbf{Proof.} By Lemma \ref{lemma1}, we have $\mathcal{T}_R(x^*)\in {\rm Ker}(\frac{\partial^2 V(x^*)}{\partial x^2})$. Then we show $\mathcal{T}_v(x^*)\in {\rm Ker}(\frac{\partial^2 V(x^*)}{\partial x^2})$. As $V(x)$ is translation invariant, for any $s\in \mathbb{R}^2$ it follows that
\begin{equation*}
    \begin{aligned}
        &\frac{\partial V(x^*+\mathbbm{1}_{m+n}\otimes s)}{\partial x}=\frac{\partial V(x^*)}{\partial x}=0,\\
        &\frac{\partial^2 V(x^*+\mathbbm{1}_{m+n}\otimes s)}{\partial x^2}=\frac{\partial^2 V(x^*)}{\partial x^2} \ge 0.
    \end{aligned}    
\end{equation*}
By Taylor expansion, we have
\begin{equation} \label{taylor}
    \begin{aligned}
        V(x^*&+\mathbbm{1}_{m+n}\otimes \nu s)=V(x^*) + \frac{\partial V(x^*)}{\partial x}(\mathbbm{1}_{m+n}\otimes \nu s) \\
        &+\nu^2(\mathbbm{1}_{m+n}\otimes s)^{\rm T}\frac{\partial^2 V(x^*)}{\partial x^2}(\mathbbm{1}_{m+n}\otimes s)+o(\nu^3).
    \end{aligned}
\end{equation}
Since $V(x)$ is translation invariant, $\frac{\partial V(x^*)}{\partial x}=0$, and the above equation holds for any $\nu$, it must follows that
\begin{equation*}
    (\mathbbm{1}_{m+n}\otimes s)^{\rm T}\frac{\partial^2 V(x^*)}{\partial x^2}(\mathbbm{1}_{m+n}\otimes s)=0
\end{equation*}
for any $s\in \mathbb{R}^2$, which implies $\mathcal{T}_v(x^*)\in {\rm Ker}(\frac{\partial^2 V(x^*)}{\partial x^2})$. It together with the fact that $\mathcal{T}_R(x^*)\in {\rm Ker}(\frac{\partial^2 V(x^*)}{\partial x^2})$ implies $\mathcal{T}_{\rm T}(x^*)\in {\rm Ker}(\frac{\partial^2 V(x^*)}{\partial x^2})$. Since ${\rm dim}(\mathcal{M}_{\rm T}(x^*))$ has the same dimension as ${\rm dim}(\mathcal{T}_{\rm T}(x^*))$, we have the result
\begin{equation} \label{mm2}
    {\rm Ker}(\frac{\partial^2 V(x^*)}{\partial x^2})=\mathcal{T}_{\rm T}(x^*).
\end{equation}

Since $\frac{\partial V(x^*)}{\partial x}=0$, and $\Delta x^{\rm T}\frac{\partial^2 V(x^*)}{\partial x^2}\Delta x>0$ for any $\Delta x\notin \mathcal{T}_{\rm T}(x^*)$, it follows from the Taylor expansion (\ref{taylor}) that $V(x)$ strictly increases along $\Delta x\notin \mathbb{R}^{2(m+n)}\backslash \mathcal{T}_{\rm T}(x^*)$. As $V(x)$ is constant on $\mathcal{M}_{\rm T}(x^*)$, we conclude that $\mathcal{M}_{\rm T}(x^*)$ is a strict local minimum manifold.

Since ${\rm Ker}(\frac{\partial^2 V(x^*)}{\partial x^2})=\mathcal{T}_{\rm T}(x^*) $ and $\mathcal{T}_v(x^*)\in {\rm Ker}(\frac{\partial^2 V(x^*)}{\partial x^2})$, by (\ref{mm2}) and the definition of $\mathcal{T}_{\rm T}(x^*) $ in (\ref{mm3}), we have
\begin{equation*}
    \frac{\partial^2 V(x^*)}{\partial x^2}(I_{m+n}\otimes L)x^*=0,
\end{equation*}
which corresponds to the fact that $\mathcal{T}_R(x^*)\in {\rm Ker}(\frac{\partial^2 V(x^*)}{\partial x^2})$.
Taking (\ref{mm4}) and $x^*=[x_a^{* \rm T}~x_p^{* \rm T}]^{\rm T}$ into the above equation gives
\begin{equation*}
    A_1(x_{a}^{*},x_{p}^{*})(I_m\otimes L)x_a^* + A_2(x_{a}^{*},x_{p}^{*})(I_n\otimes L)x_p^*=0 \\
\end{equation*}
It is evident that $A_1(x_{a}^{*},x_{p}^{*})(I_m\otimes L)x_a^*=0$ if and only if (\ref{mm1}) is satisfied. Let $\mathcal{T}_R(x_a^*)=\{(I_m\otimes L)x_a^*:L \in \mathfrak{so}(2)\}$, which is the tangent space of the manifold $\mathcal{M}_R(x_a^*)$, then we have ${\rm Ker}(A_1(x_a^*,x_p^*))=\mathcal{T}_R(x_a^*)$. Therefore, the number of zero eigenvalues of $A_1(x_a^*,x_p^*)$ equals the dimension of $\mathcal{M}_R(x_a^*)$. Similarly, $A_3(x_{a}^{*},x_{p}^{*})(I_n\otimes L)x_p^*=0$ if and only if (\ref{mm6}) is satisfied, then, the number of zero eigenvalues of $A_3(x_a^*,x_p^*)$ equals ${\rm dim}(\mathcal{M}_R(x_p^*))$.
$\hfill\square$

\begin{lemma} \label{lemma2}
    Consider the symmetric circulant matrix with rotation 
    \begin{equation*}
    H=\left[\begin{array}{ccccc}
    C_{0} & C_{1} & \cdots & C_{m-1} \\
    R C_{m-1} R^{-1} & R C_{0} R^{-1} & \cdots & R C_{m-2} R^{-1} \\
    R^{2} C_{m-2} R^{-2} & R^{2} C_{m-1} R^{-2}  & \cdots & R^{2} C_{m-3} R^{-2} \\
    & \vdots & & \\
    R^{m-1} C_{1} R^{1-m} &  R^{m-1} C_{2} R^{1-m} &    \cdots & R^{m-1} C_{0} R^{1-m}
    \end{array}\right]
    \end{equation*}
where $R\in SO(2)$, $R^m=I$, and each block matrix $C_i\in \mathbb{R}^{2\times 2}$, $i=0,\cdots,m-1$, is also symmetric. Let
\begin{equation*}
    D_k=C_0+w^kR^{m-1}C_1 + w^{2k}R^{m-2}C_2+\cdots+ w^{(m-1)k}RC_{m-1}
\end{equation*}
where $k=0,1,2,\cdots,m-1$, and $w=e^{-2\pi {\rm i }/m}$ is the primitive $n$th root of unity. Then, the eigenvalues of $H$ consist of the eigenvalues of $D_k$, and the eigenvectors are 
\begin{equation*}
    v_{k,l}\!=\!\begin{bmatrix}
        \psi_l^{\rm T} & w^{k(m\!-\!1)}\psi_l^{\rm T}R^{\rm T} &  w^{k(m\!-\!2)} \psi_l^{\rm T}R^{2 \rm T}& \cdots & w^k \psi_l^{\rm T}R^{m\!-\!1,\rm T}
    \end{bmatrix}^{\rm T}
\end{equation*}
where $\psi_l$, $l=1,2$ is the eigenvector of $D_k$.    
\end{lemma}
The proof can be found in appendix.

\subsection{Topology design}\label{sec3.2}

The inter-agent topologies for the active and passive groups are designed as follows. 
\newtheorem{assumption}{Assumption}
\begin{assumption} \label{ass1}
(Topology of the active group)
    \begin{itemize}
    \item Case (i). If $m$ is odd, the inter-agent topology is given by an undirected ring. In this case, we have $d_i^a=2$ for all $i\in \mathcal{V}^a$;
    \item Case (ii). If $m=2l$ and $l$ is odd, the agents are divided into two distinct groups, each containing $l$ players. The agents in each subgroup are connected by a ring. Besides, the two subgroups are linked by $l$ one-to-one edges, each connecting a distinct pair of agents from the two subgroups.
    In this case, $d_i^a=3$ for all $i\in \mathcal{V}^a$;
    \item Case (iii). If $m=4l$ and $l=1,2,\cdots$, we divide the vertices into $l$ distinct subgroups with four vertices each. The connectivity of each subgroup is specified in Algorithm~\ref{alg1} where $GP_i:=\{gp_{i,1},gp_{i,2},gp_{i,3},gp_{i,4}\}$, $i=1,2,\cdots,l$ denotes the set of players in the $i$th subgroup. Note that agents in each subgroup are completely connected, and each agent connects to two agents from every other subgroup. In this case, $d_i^a=1+m/2$ for all $i\in \mathcal{V}^a$.
\end{itemize}
\end{assumption}

\begin{algorithm}
    \caption{Topology design for the case $m=4l$, $l=1,2,\cdots$}
\begin{algorithmic}[1] \label{alg1}
    \STATE {\textbf{Input:} $l$, $GP_i$ for $i=1,2,\cdots,l$}
    \IF{$l=1$}
        \STATE {the topology is a complete graph}
    \ELSE
        \FOR{$i=2:l$}
            \STATE{the topology of subgroup $i$ is a complete graph}
            \FOR{$k=1:i-1$}
                \STATE{player $gp_{k,1}$ connects with $gp_{i,2}$ and $gp_{i,3}$}
                \STATE{player $gp_{k,2}$ connects with $gp_{i,3}$ and $gp_{i,4}$}
                \STATE{player $gp_{k,3}$ connects with $gp_{i,1}$ and $gp_{i,4}$}
                \STATE{player $gp_{k,4}$ connects with $gp_{i,1}$ and $gp_{i,2}$}
            \ENDFOR
        \ENDFOR
    \ENDIF
\end{algorithmic}
\end{algorithm}

\begin{assumption}\label{ass2}
(Topology of the passive group)
For the passive group, the inter-agent topology is given by an undirected ring, regardless of whether $n$ is odd or even.
\end{assumption}

Here, we provide some insights into the design of the inter-agent topologies. Without a repelling term within the active group, its agents would converge to the same point as they approach the passive group, causing the encirclement control to fail. Meanwhile, the passive group may diverge without attractive forces from its neighbors, but only two neighbors are sufficient to maintain cohesion. This difference is the reason for the distinct topologies designed for the two groups. When $m$ is odd, an undirected ring topology is sufficient because players cannot be paired to reach consensus due to the symmetry of the configuration. This also explains why, in the case where $m=2l$ and $l$ is odd, the graph is a ring in each subgroup. We connect $l$ distinct pairs of players from the two subgroups to avoid agent overlap. The case (iii) is more complex as players are more likely to converge to the same point. The condition that each subgroup is completely connected ensures players in each sub-group are evenly distributed on the circle, and the interactions among subgroups prevents some players from reaching the same point.

In the following Section~\ref{sec4} we will see that the convergence property of the game is closely related to the positive definiteness of the following $\mathcal{R}$-invariant and translation invariant potential functions
\begin{align}
    W_i^a(x) &\!=\! \frac{\alpha_1}{2n}\sum_{i=1}^m\sum_{j=1}^n \!\|x_i^a \!-\! x_j^p\|^2 \!+\! \frac{\alpha_2}{d_i^a}\sum_{i=1}^m\sum_{k\in \mathcal{N}_i^a} \!\frac{1}{\|x_i^a\!-\!x_k^a\|}  \label{Wpx}\\
    W_j^p(x) &\!=\! \beta_1\sum_{j=1}^n \! \frac{1}{\|x_j^p\!-\!\frac{1}{m}\sum_{i=1}^m \! x_i^a\|} \!+\! \frac{\beta_2}{6}\sum_{j=1}^n \!\sum_{k\in \mathcal{N}_j^p} \!\!\|x_j^p\!-\!x_k^p\|^3  \label{Wex}
\end{align}
where $\alpha_1$, $\alpha_2$, $\beta_1$ and $\beta_2$ are positive constants, $d_i^a$ is the number of active neighbors of active player $i$,
$x_i^a\ne x_k^a$ for any $i\ne k,~i,k\in \mathcal{V}^a$, and $x_j^p\ne \frac{1}{m}\sum_{i=1}^m x_i^a$ for each $j\in \mathcal{V}^p$.

\begin{lemma} \label{lemma3}
    Under Assumptions~\ref{ass1} and \ref{ass2}, for any positive coefficients $\alpha_1$, $\alpha_2$, $\beta_1$ and $\beta_2$, $W_i^a(x)$ and $W_j^p(x)$ are positive semi-definite at two concentric circles with radii given in Table~\ref{table1}, the Hessian of $W_i^a(x)$ and $W_j^p(x)$ each have exactly 3 zero-eigenvalues, and both functions obtain local minimum at $\mathcal{M}_{\rm T}(x^*)$. Furthermore, 
    the players are distributed on the circles as follows.
    For the active group, the active players are evenly located on one circle, and for any active player $i$, its neighbors are positioned at the farthest points from it; for the passive group, the passive players are evenly located on the other circle, and for any passive player $j$, its neighbors are positioned at the nearest points from it.
\end{lemma}

\begin{table}[h] 
	\begin{center}
		\caption{Radii of circles}
	\begin{tabular}{|c|c|c|}
        \hline
        ~& $r_a$ & $r_p$ \\
        \hline 
        $m$ is odd & $\sqrt[3]{\frac{\alpha_2(1-cos\frac{(m-1)\pi}{m})}{4\alpha_1sin\frac{(m-1)\pi}{2m}}}$ & \multirow{3}{*}{$\sqrt[4]{\frac{\beta_1}{4\beta_2 sin\frac{\pi}{n}(1-cos\frac{2\pi}{n})}}$}  \rule{0pt}{17pt}\\
        \cline{1-2}
        \makecell{$m=2l$, \\$l$ is odd}        & $\sqrt[3]{\frac{\alpha_2(1+cos^3\frac{\pi}{m}+cos\frac{2\pi}{m})}{4\alpha_1cos^3\frac{\pi}{m}}}$ &~ \rule{0pt}{17pt}\\
        \cline{1-2}
        $m=4l$ & $\sqrt[3]{\frac{\alpha_2}{\alpha_1}(\frac{1}{4}+\sum_{l=\frac{m}{4}}^{\frac{m}{2}-1}\frac{1-cos\frac{2\pi l}{m}}{4 sin^3 \frac{\pi l}{m}})}$ & ~ \rule{0pt}{17pt}\\
        \hline  
    \end{tabular}
    \label{table1}
    \end{center}
\end{table}

\textbf{Proof.} By Theorem \ref{them1}, the claim can be proved by just checking one arbitrary point on the manifold. We first verify the properties related to $W_i^a(x)$ in detail, then similar method can be used to show that of $W_j^p(x)$.

Without loss of generality, we choose a coordinate of vertices 
\begin{equation*}
\begin{aligned}
 x_i^{a*}&=[r_acos(2(i-1)\pi/m)~~r_asin(2(i-1)\pi/m)]^{\rm T}, ~i\in \mathcal{V}^a   \\
 x_j^{p*}&=[r_pcos(2(j-1)\pi/n)~~r_psin(2(j-1)\pi/n)]^{\rm T},~j\in \mathcal{V}^p. 
\end{aligned}
\end{equation*}

\begin{figure}	
	\centering
	\includegraphics[width=0.14\textwidth]{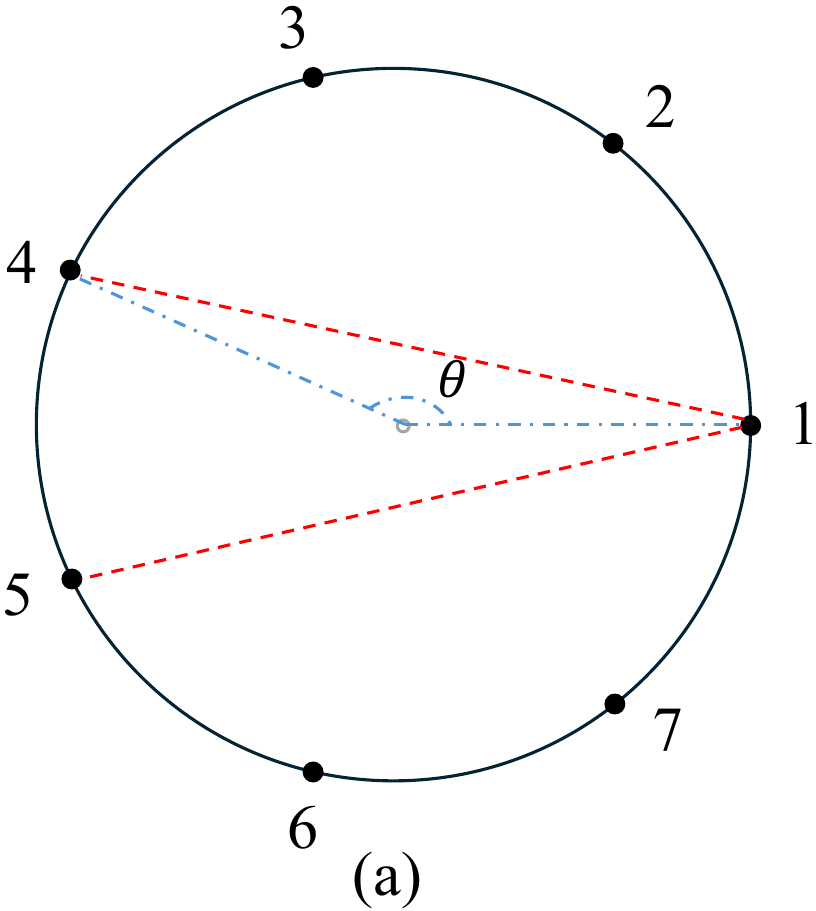}~~
        \includegraphics[width=0.14\textwidth]{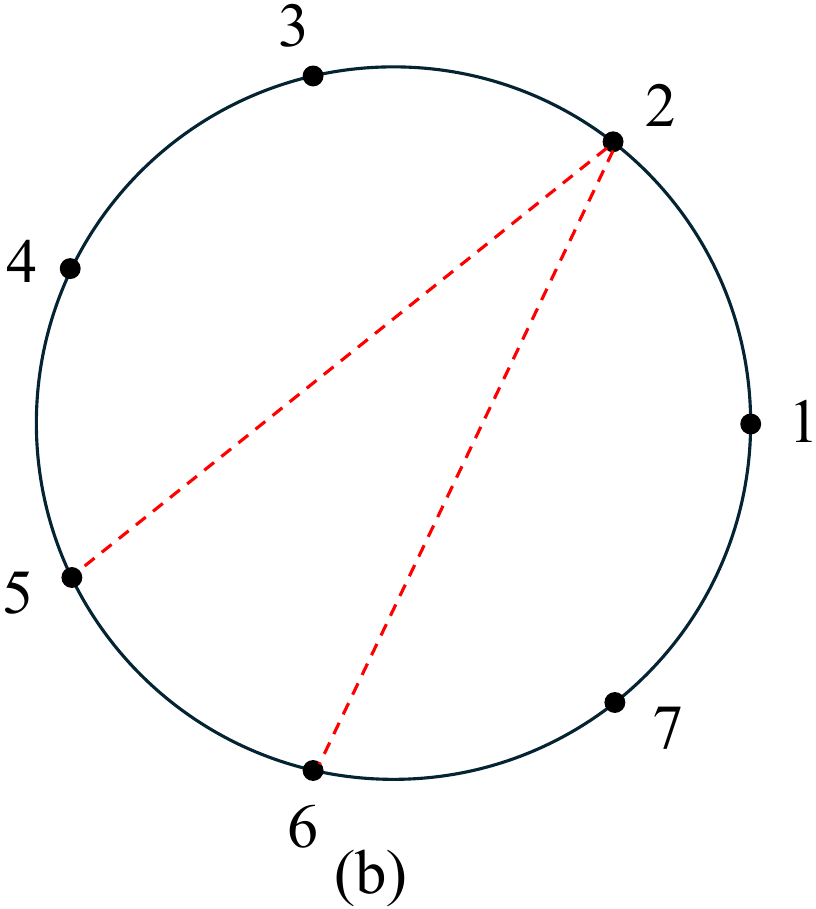}~~
        \includegraphics[width=0.14\textwidth]{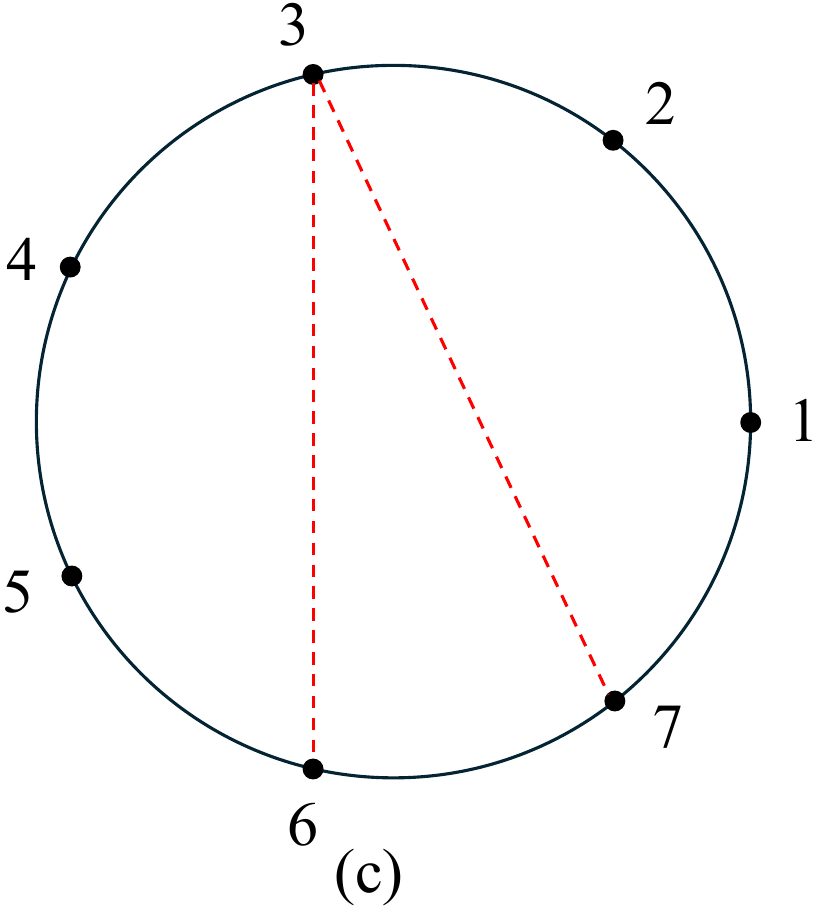}\\
        \includegraphics[width=0.19\textwidth]{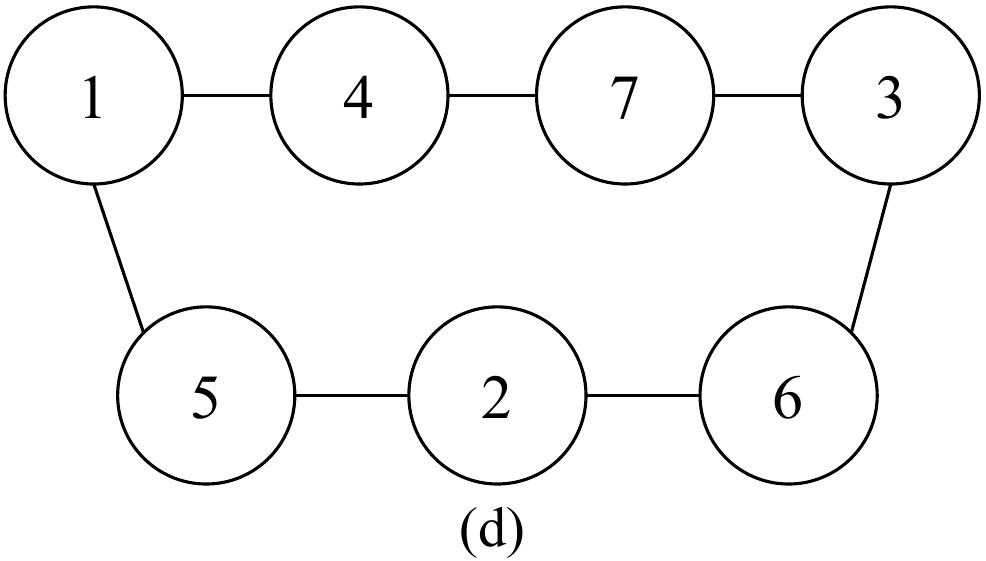}
	\caption{Illustration of the relative positions of active players under the topology (d), when $m=7$.}
	\label{circles}
\end{figure}

Denote $x_i^a=[x_{i1}^a~x_{i2}^a]^{\rm T}$ and $x_j^p=[x_{j1}^p~x_{j2}^p]^{\rm T}$. Direct computation gives the Hessian of $W_i^a(x)$ by
\begin{equation}
            \frac{\partial^{2} W_i^a(x^{*})}{\partial x^{2}}=\left[\begin{array}{c|c}
    G_1^a & G_2^a \\
    \hline
    G_2^{a \rm T} & G_3^a 
    \end{array}\right]
\end{equation}
where $G_1^a\in\mathbb{R}^{2m\times 2m}$, $G_2^a\in \mathcal{R}^{2m\times 2n}$, $G_3^a\in \mathcal{R}^{2n\times 2n}$, with $G_{1,ik}^a$ following (\ref{mm7}) for $l,k=1,2,\cdots,m$,  
\begin{figure*}
    \begin{equation} \label{mm7}
        G_{1,lk}^a=\left\{
        \begin{aligned}
        &\begin{bmatrix}
            \alpha_1+\alpha_2\sum_{q\in \mathcal{N}_l^a}\frac{2(x_{l1}^{a*}-x_{q1}^{a*})^2-(x_{l2}^{a*}-x_{q2}^{a*})^2}{\|x_l^{a*}-x_q^{a*}\|^5} & -\alpha_2 \sum_{q\in \mathcal{N}_l^a}\frac{3(x_{l1}^{a*}-x_{q1}^{a*})(x_{l2}^{a*}-x_{q2}^{a*})}{\|x_l^{a*}-x_q^{a*}\|^5} \\
            -\alpha_2 \sum_{q\in \mathcal{N}_l^a}\frac{3(x_{l1}^{a*}-x_{q1}^{a*})(x_{l2}^{a*}-x_{q2}^{a*})}{\|x_l^{a*}-x_q^{a*}\|^5} & \alpha_1+\alpha_2\sum_{q\in \mathcal{N}_l^a}\frac{2(x_{l2}^{a*}-x_{q2}^{a*})^2-(x_{l1}^{a*}-x_{q1}^{a*})^2}{\|x_l^{a*}-x_q^{a*}\|^5}
        \end{bmatrix}, &{\rm if} ~~k=l, \\
        &\begin{bmatrix}
            -\alpha_2 \frac{2(x_{l1}^{a*}-x_{k1}^{a*})^2-(x_{l2}^{a*}-x_{k2}^{a*})^2}{\|x_l^{a*}-x_k^{a*}\|^5} & -\alpha_2\frac{3(x_{l1}^{a*}-x_{k1}^{a*})(x_{l2}^{a*}-x_{k2}^{a*})}{\|x_l^{a*}-x_k^{a*}\|^5} \\
            -\alpha_2 \frac{3(x_{l1}^{a*}-x_{k1}^{a*})(x_{l2}^{a*}-x_{k2}^{a*})}{\|x_l^{a*}-x_k^{a*}\|^5} & -\alpha_2 \frac{2(x_{l2}^{a*}-x_{k2}^{a*})^2-(x_{l1}^{a*}-x_{k1}^{a*})^2}{\|x_l^{a*}-x_k^{a*}\|^5}
        \end{bmatrix}, &{\rm if} ~~k\in \mathcal{N}_l^a, \\
        & 0_{2\times 2}, &{\rm otherwise},
       \end{aligned} 
        \right.
    \end{equation}
\end{figure*}

\begin{equation*}
\begin{aligned}
    G_{2,lk}^a&=-\frac{\alpha_1}{n}I_2, ~l=1,2,\cdots,m,~k=1,2,\cdots,n,\\
    G_3^a&=\frac{m}{n}\alpha_1I_{2n}.
\end{aligned}
\end{equation*}

Then $\frac{\partial^{2} W_i^a(x^{*})}{\partial x^{2}}$ has the following properties. (1) Since the topology is undirected, $\frac{\partial^{2} W_i^a(x^{*})}{\partial x^{2}}$ is symmetric, and each block matrix of $G_1^a,~G_2^a,~ G_3^a$ is also symmetric. This can also be verified from the explicit formulas given above. (2) The sum of each row is zero. (3) As the passive players are evenly distributed on the circle, $\sum_{i=1}^m x_{i1}^{p*}=0$, $\sum_{i=1}^m x_{i2}^{p*}=0$, which leads to $G_2^a(I_n\otimes L)x_p^*=0$.
Therefore, by Theorem~\ref{them1}, to show $W_i^a(x)$ obtains a strict local minimum at $\mathcal{M}_{\rm T}(x^*)$, i.e., the number of zero eigenvalues of $\frac{\partial^{2} W_i^a(x^{*})}{\partial x^{2}}$ equals 3, we only need to show $G_1^a$ is positive semi-definite and has exactly one zero eigenvalue. 

To this end, we first show that $G_1^a$ is positive semi-definite. Without loss of generality, we take Case (i) of Assumption~\ref{ass1} that $m$ is odd as an example, where $x_1^{a*}=[r_a~0]^{\rm T}$, and the neighbors of active player 1 are active players $(m+1)/2$ and $(m+3)/2$, as they are located at the farthest point. 
For example, when $m=7$ in Fig.~\ref{circles}, the neighbors of active player 1 are players 4 and 5, which positioned at the farthest points on the circle to player 1.
Denote $\theta=(m-1)\pi/m$, then $x_{(m+1)/2}^{a*}=[r_acos\theta~r_asin\theta]^{\rm T}$ and $x_{(m+3)/2}^{a*}=[r_acos\frac{(m+1)\pi}{m}~r_asin\frac{(m+1)\pi}{m}]^{\rm T}$. By Table~\ref{table1}, the radius is $r_a=\sqrt[3]{\frac{\alpha_2(1-cos\theta)}{4\alpha_1 sin \frac{\theta}{2}}}$.
For simplicity, let $x_k^{a*}=[r_acos\theta_k~r_asin\theta_k]^{\rm T}$ where $k=(m+1)/2$ or $k=(m+3)/2$. A direct calculation gives $2(x_{i1}^{a*}-x_{k1}^{a*})^2-(x_{i2}^{a*}-x_{k2}^{a*})^2=(1-cos\theta_k)(1-3cos\theta_k)$ which is positive if and only if $cos\theta_k<1/3$. This condition is satisfied when $\theta_k=(m-1)\pi/m$ or $(m+1)\pi/m$. 
Given $x_k^{a*}$ and $r_a$, we have $\alpha_1+\alpha_2\sum_{k\in \mathcal{N}_i^a}\frac{2(x_{12}^{a*}-x_{k2}^{a*})^2-(x_{11}^{a*}-x_{k1}^{a*})^2}{\|x_1^{a*}-x_k^{a*}\|^5}=\alpha_1+\alpha_1\frac{sin^3 \theta(1+3cos\theta)}{4sin^5 (\theta/2)}$, which is also positive. Thus, all the diagonal elements of $G_1^a$ are positive. This fact, together with the observation that each row sum of $G_1^a$ is $\alpha_1$, leads to the conclusion that $G_1^a$ is positive semi-definite.

Now, we show that $G_1^a$ has one zero eigenvalue.
It is noted that $G_1^a$ is independent of $x_p^*$. As shown in Fig.~\ref{circles}, the subgroup $\{2,5,6\}$ can be obtained by rotating $\{1,4,5\}$, and further rotation leads to $\{3,6,7\}$.
By Definition 1, since $W_i^a(x)$ is $\mathcal{R}$-invariant, it follows that $\frac{\partial^2 W_i^a(x)}{\partial x^2}|_{x=\bar Rx}=\bar R \frac{\partial^2 W_i^a(x)}{\partial x^2} \bar R^{\rm T}$. Therefore, $G_1^a$ can be written as 
\begin{equation*}
        G_1^a\!=\!\!\begin{bmatrix}
        E_0 & E_1  & \cdots & E_{m-1} \\
        \bar R E_{m-1} \bar R^{\rm T} & \bar R E_0 \bar R^{\rm T} &  \cdots & \bar R E_{m-2} \bar R^{\rm T} \\
        \vdots & \vdots & \vdots & \vdots  \\
        \bar R^{m-1} E_1 \bar R^{1-m} & \bar R^{m-1} E_2 \bar R^{1-m}  & \cdots &
        \bar R^{m-1} E_0 \bar R^{1-m} &
    \end{bmatrix}
\end{equation*}
where $\bar R\in SO(2)$, $\bar R^m=I$. For $i=1,2,\cdots,m$, $E_i\in \mathbb{R}^{2\times 2}$ can be obtained from (\ref{mm7}). By Lemma~\ref{lemma2}, the eigenvalues of $G_1^a$ are given by the eigenvalues of 
$$D_k=E_0+w^kR^{m-1}E_1 + w^{2k}R^{m-2}E_2+\cdots+ w^{(m-1)k}RE_{m-1},$$ 
for $k=0,1,\cdots,m-1$. We first compute the eigenvalues of $D_0$, again considering Case (i) of Assumption~\ref{ass1} as an example. Here, we have the rotation matrix
$\bar R=\begin{bmatrix}
    cos(2\pi/m) & -sin(2\pi/m)\\
    sin(2\pi/m) & cos(2\pi/m)
\end{bmatrix}.$


It is straightforward to see from (\ref{mm7}) that 
$$E_0\!=\!\alpha_1 I_2\!-\!E_1\!-\!\cdots\!-\!E_{m-1}\!=\!\alpha_1 I_2\!-\!E_{(m-1)/2}\!-\!E_{(m+1)/2},$$ we thus have
\begin{equation} \label{mm8}
    D_0=\alpha_1 I_2 + (R^{\frac{m+1}{2}}-I_2)E_{\frac{m-1}{2}}+ (R^{\frac{m-1}{2}}-I_2)E_{\frac{m+1}{2}}.
\end{equation}
Let $\theta=(m-1)\pi/m$.
Then, $\|x_1^{a*}-x_{(m+1)/2}^{a*}\|=\|x_1^{a*}-x_{(m+3)/2}^{a*}\|=2r_asin\frac{\theta}{2}$. Taking this into (\ref{mm8}) gives
\begin{equation*}
\begin{aligned}
    D_0& \!=\! \alpha_1 \Big(I_2\!+\!\frac{1}{16sin^4(\frac{\theta}{2})}\big(
    \begin{bmatrix}
        4(1\!-\!cos\theta)^2 & 4sin\theta(1\!-\!cos\theta) \\
        2sin\theta(1\!-\!cos\theta) & -2(1\!-\!cos\theta)^2
    \end{bmatrix} \\
    & \qquad \qquad \quad+
    \begin{bmatrix}
        4(1-cos\theta)^2 & -4sin\theta(1-cos\theta) \\
        -2sin\theta(1-cos\theta) & -2(1-cos\theta)^2
    \end{bmatrix}
    \big) \Big) \\
    &=\alpha_1 \begin{bmatrix}
        3 & 0 \\0 & 0
    \end{bmatrix}
\end{aligned}
\end{equation*}
which has exactly one zero eigenvalue. Similarly, we can prove that the matrices $D_k$, $k=1,2,\cdots,m-1$ are nonsingular. Therefore, the number of zero eigenvalues of $G_1^a$ is exactly 1. The same conclusion holds for Cases (ii) and (iii).

For the passive group, the Hessian matrix of $W_j^p(x)$ is
\begin{equation*}
        \frac{\partial^{2} W_j^p(x^{*})}{\partial x^{2}}=\left[\begin{array}{c|c}
    G_1^p  & G_2^p \\
    \hline
    G_2^{p\rm T} & G_3^p 
    \end{array}\right]
\end{equation*}
where each block of $G_1^p$ is 
\begin{equation*}
        \begin{bmatrix}
            \!-\frac{\beta_1n}{m^2r_p^3}\!+\!\frac{3\beta_1}{m^2r_p^5} \!\sum\limits_{j=1}^n(x_{j1}^{p*}\!-\!x_{c1}^{a*})^2 & \frac{3\beta_1}{m^2r_p^5}\!\sum\limits_{j=1}^n(x_{j2}^{p*}\!-\!x_{c2}^{a*})^2 \\
            \!\frac{3\beta_1}{m^2r_p^5}\!\sum\limits_{j=1}^n(x_{j2}^{p*}\!-\!x_{c2}^{a*})^2 & -\frac{\beta_1n}{m^2r_p^3}\!+\!\frac{3\beta_1}{m^2r_p^5}\!\sum\limits_{j=1}^n(x_{j2}^{p*}\!-\!x_{c2}^{a*})^2
        \end{bmatrix}\!,
\end{equation*}
the elements of $G_2^p$ is, for $l=1,\cdots,m$ and $k=1,\cdots,n$,
\begin{equation*}
G_{2,lk}^p\!=\!\!
    \begin{bmatrix}
        \frac{\beta_1}{m}(\!-\!\frac{3}{r_p^5}(x_{k1}^{p*}\!-\!x_{c1}^{a*})^2\!+\!\frac{1}{r_p^3}) & -\!\frac{\beta_1}{m}\frac{3}{r_p^5}(x_{k2}^{p*}\!-\!x_{c2}^{a*})^2 \\
        -\!\frac{\beta_1}{m}\frac{3}{r_p^5}(x_{k2}^{p*}\!-\!x_{c2}^{a*})^2 & \frac{\beta_1}{m}(-\!\frac{3}{r_p^5}(x_{k2}^{p*}\!-\!x_{c2}^{a*})^2\!+\!\frac{1}{r_p^3})
    \end{bmatrix}\!,
\end{equation*}
and $G_3^p$ is given by (\ref{mm9}), where $x_c^a=[x_{c1}^a~x_{c2}^a]^{\rm T}:=\frac{1}{m}\sum_{i=1}^m x_i^a$ is the center of the active group. We can verify that $\frac{\partial^2 W_j^p(x)}{\partial x^2}$ also satisfies the properties of  $\frac{\partial^2 W_i^a(x)}{\partial x^2}$ under Assumption~\ref{ass2}. (1) It is symmetric. (2) The sum of each row is zero. (3) On the desired manifold, 
$G_2^{p\rm T}(I_m\otimes L)x_a^*=0$. (4) On the desired manifold where neighbors of each passive player are located at the nearest points from that player, $G_3^p\ge 0$ and has one zero eigenvalue. Therefore, by Theorem~\ref{them1}, $W_j^p(x)$ obtains local minimum at $\mathcal{M}_{\rm T}(x^*)$. 
\begin{figure*}
    \begin{equation} \label{mm9}
        G_{3,lk}^p\!=\!\left\{
        \begin{aligned}
        &\begin{bmatrix}
            \!\beta_1 (\frac{3}{r_p^5}(x_{l1}^{p*}\!-\!x_{c1}^{a*})^2\!-\!\frac{1}{r_p^3}) \!+\! \beta_2 \!\!\sum\limits_{q\in\mathcal{N}_l^p} \!(\frac{(x_{l1}^{p*}\!-\!x_{q1}^{p*})^2}{\|x_l^{p*}\!-\!x_q^{p*}\|}\!+\!\|x_l^{p*}-x_q^{p*}\|) & \frac{3\beta_1 }{r_p^5}(x_{l1}^{p*}\!-\!x_{c1}^{a*})(x_{l2}^{p*}\!-\!x_{c2}^{a*})\!+\!\beta_2 \!\!\sum\limits_{q\in\mathcal{N}_l^p}\!\!\frac{(x_{l1}^{p*}\!-\!x_{q1}^{p*})(x_{l2}^{p*}\!-\!x_{q2}^{p*})}{\|x_l^{p*}\!-\!x_q^{p*}\|} \\
            \frac{3\beta_1 }{r_p^5}(x_{l1}^{p*}\!-\!x_{c1}^{a*})(x_{l2}^{p*}\!-\!x_{c2}^{a*})\!+\!\beta_2 \!\!\sum\limits_{q\in\mathcal{N}_l^p}\!\!\frac{(x_{l1}^{p*}\!-\!x_{q1}^{p*})(x_{l2}^{p*}\!-\!x_{q2}^{p*})}{\|x_k^{p*}\!-\!x_q^{p*}\|} & \beta_1 (\frac{3}{r_p^5}(x_{l2}^{p*}\!-\!x_{c2}^{a*})^2\!-\!\frac{1}{r_p^3}) \!+\! \beta_2\!\!\sum\limits_{q\in\mathcal{N}_l^p}\!\!(\frac{(x_{l2}^{p*}\!-\!x_{q2}^{p*})^2}{\|x_l^{p*}\!-\!x_q^{p*}\|}\!+\!\|x_l^{p*}\!-\!x_q^{p*}\|)
        \end{bmatrix}, \\
        &\qquad \qquad \qquad \qquad \qquad \qquad \qquad \qquad \qquad \qquad \qquad \qquad \qquad \qquad \quad~~  {\rm if} ~k=l, \\
        &\begin{bmatrix}
            -\beta_2(\frac{(x_{l1}^{p*}-x_{k1}^{p*})^2}{\|x_l^{p*}-x_k^{p*}\|}+\|x_l^{p*}-x_k^{p*}\|) & -\beta_2 \frac{(x_{l1}^{p*}-x_{k1}^{p*})(x_{l2}^{p*}-x_{k2}^{p*})}{\|x_l^{p*}-x_k^{p*}\|} \\
            -\beta_2 \frac{(x_{l1}^{p*}-x_{k1}^{p*})(x_{l2}^{p*}-x_{k2}^{p*})}{\|x_l^{p*}-x_k^{p*}\|} & -\beta_2(\frac{(x_{l2}^{p*}-x_{k2}^{p*})^2}{\|x_l^{p*}-x_k^{p*}\|}+\|x_l^{p*}-x_k^{p*}\|)
        \end{bmatrix}, ~~~{\rm if} ~~k\in \mathcal{N}_l^p, \\
        & 0_{2\times 2}, \qquad\qquad\qquad\qquad\qquad\qquad\qquad\qquad\qquad\qquad\qquad\qquad\quad ~\qquad {\rm otherwise}.
       \end{aligned} 
        \right.
    \end{equation}
\end{figure*}
$\hfill\square$

Now we can present an alternative approach for designing the graphs, given the desired ordering of players on the circle. By Lemma~\ref{lemma3}, the neighbors of each active player are positioned at the farthest points from it, with the number of neighbors specified in Assumption~\ref{ass1}. Therefore, given the desired ordering of active players on the circle, each player connects with its $d_i^a$ farthest players. Similarly, given the desired ordering of passive players on the circle, the neighbors of each passive player are the nearest two players.

\section{Game analysis} \label{sec4}
In this section, the encirclement problem is formulated as an infinite time-horizon differential game and asymptotic convergence to the desired manifold is proved. 

For each individual cost in (\ref{obj_a}) and (\ref{obj_p}) we choose $q_i^a(x)$ and $q_j^p(x)$ as a trade-off between chasing or escaping and inter-agent interactions. Specifically, for $i\in\mathcal{V}^a$ and $j\in\mathcal{V}^p$, we define
    \begin{align}
        q_i^a(x)&\!=\!x^{\rm T}C^{a{\rm T}}(x)Q_i^aC^a(x)x \!+\! 2\!\!\!\sum_{j=m+1}^{m+n} \!\!\!x^{\rm T}C^{a{\rm T}}(x)B_jB_j^{\rm T}C^p(x) x \label{mm10}\\
        q_j^p(x)&\!=\!x^{\rm T}C^{p{\rm T}}(x)Q_j^pC^p(x)x \!+\! 2\!\sum_{i=1}^{m} \!x^{\rm T}C^{p{\rm T}}(x)B_iB_i^{\rm T}C^a(x) x \label{mm16}
    \end{align}
where $Q_i^a\in \mathbb{S}_+^{2(m+n)\times 2(m+n)}$ and $Q_j^p\in \mathbb{S}_+^{2(m+n)\times 2(m+n)}$ are some constant matrices, and the block matrices $C^a(x)$ and $C^p(x)$ satisfy
\begin{equation} \label{mm11}
    C^a(x)=\left[\begin{array}{c|c}
    C_1^a(x) & C_2^a(x) \\
    \hline
    C_2^{a {\rm T}}(x) & C_3^a(x) 
    \end{array}\right],~~
    C^p(x)=\left[\begin{array}{c|c}
    C_1^p(x) & C_2^p(x) \\
    \hline
    C_2^{p {\rm T}}(x) & C_3^p(x) 
    \end{array}\right]
\end{equation}
with 
\begin{equation*}
\begin{aligned}
    C_{1,lk}^a&=\left\{
    \begin{aligned}
      &(\alpha_1-\alpha_2\sum_{q\in \mathcal{N}_l^a}\frac{1}{\|x_l^a-x_q^a\|^3})I_2 & {\rm if}~k=l\\
      &\alpha_2\frac{1}{\|x_l^a-x_k^a\|^3}I_2 & {\rm if}~ k\in \mathcal{N}_l^a \\
      & 0_{2\times 2} & {\rm otherwise}
    \end{aligned}  
    \right.\\
   C_{2,lk}^a&= -\frac{\alpha_1}{n}I_2,~~ l=1,\cdots,m,~k=1,\cdots,n, \\
   C_{3,lk}^a&= \frac{m}{n}\alpha_1 I_{2n},~~l,k=1,\cdots,n.
\end{aligned}
\end{equation*}
where $\alpha_1>0$ and $\alpha_2>0$ are constants.

The intuition behind the objective function $q_i^a$ of the active group is to drive its members to move towards the passive group while maintaining sufficient separation between them. This ensures that, at the desired configuration, they face balanced forces so the configuration does not change. At the same time, the active players aim to draw the passive players closer to them. This is why $q_i^a$ consists of two parts.

Similarly, the elements of the block matrix $C^p(x)$ are
\begin{equation*}
\begin{aligned}
    C_{1,lk}^p&=-\frac{\beta_1}{m^2}\sum_{l=1}^m\frac{1}{\|x_l^p-x_c^a\|^3}I_2,~l,k=1,\cdots,m\\
    C_{2}^p&=\frac{\beta_1}{m} \!\begin{bmatrix}
        \frac{1}{\|x_1^p-x_c^a\|^3}\mathbbm{1}_m & \frac{1}{\|x_2^p-x_c^a\|^3}\mathbbm{1}_m &  \cdots
        &  \frac{1}{\|x_n^p-x_c^a\|^3}\mathbbm{1}_m
    \end{bmatrix}\!\otimes\! I_2\\
    C_{3,lk}^p&\!=\! \left\{ \!
    \begin{aligned}
        &(-\beta_1\frac{1}{\|x_l^p\!-\!x_c^a\|^3} + \beta_2 \!\sum_{q\in \mathcal{N}_l^p}\!\|x_l^p\!-\!x_q^p\|)I_2,  ~~  {\rm if}~k=l\\
        &(-\beta_2\sum_{k\in \mathcal{N}_l^p}\|x_l^p-x_k^p\|)I_2, \qquad \qquad \qquad  {\rm if} ~k\in \mathcal{N}_l^p
    \end{aligned}
    \right.
\end{aligned}
\end{equation*}
where $x_c^a:=\frac{1}{m}\sum_{i=1}^m x_i^a$ is the center of the active group.

The intuition behind the objective function $q_j^p$ of the passive group is to make its players flee from the active group while attracting their neighbors to maintain sufficient separation. This ensures that they converge to the desired configuration, in which they also face balanced forces. Simultaneously, the passive players drive the active players away from them. Thus, $q_j^p$ also consists of two terms. 

In order to obtain the Nash equilibrium strategies for the differential games (\ref{game_a}) and (\ref{game_p}), the following coupled Hamilton-Jacobi-Bellman (HJB) equations should be satisfied:
    \begin{align}
         \frac{1}{2}q_i^a(x)-\frac{1}{2}\frac{\partial V_i^{a {\rm T}}}{\partial x} &B_iB_i^{\rm T}\frac{\partial V_i^{a}}{\partial x} 
        -\sum_{k=1,k\ne i}^m\frac{\partial V_i^{a {\rm T}}}{\partial x} B_kB_k^{\rm T}\frac{\partial V_k^{a}}{\partial x}  \notag\\
         &-\sum_{j=1}^n \frac{\partial V_i^{a {\rm T}}}{\partial x} B_{j+m}B_{j+m}^{\rm T}\frac{\partial V_j^{p}}{\partial x} =0, \label{HJBi}\\
         \frac{1}{2}q_j^p(x)-\frac{1}{2}\frac{\partial V_j^{p {\rm T}}}{\partial x} &B_{j+m}B_{j+m}^{\rm T}\frac{\partial V_j^{p}}{\partial x} 
        \!-\!\!\!\sum_{\substack{k=m+1,\\k\ne m+j}}^{m+n} \!\!\!\frac{\partial V_j^{p {\rm T}}}{\partial x} B_{k}B_{k}^{\rm T}\!\frac{\partial V_k^{p}}{\partial x} \notag\\
         &-\sum_{i=1}^m \frac{\partial V_j^{p {\rm T}}}{\partial x} B_iB_i^{\rm T}\frac{\partial V_i^{a}}{\partial x} =0, \label{HJBj}
    \end{align}
with $i=1,\cdots,m$, $j=1,\cdots,n$, and $V_i^a(x)$ and $V_j^p(x)$ are locally positive definite at $\mathcal{M}_{\rm T}(x)$.

Provided that solutions $V_i^a(x)$ and $V_j^p(x)$ to the equations (\ref{HJBi}) and (\ref{HJBj}) exist, the Nash equilibrium strategies are given by
\begin{equation} \label{nash}
    \begin{aligned}
        u_i^a&=-B_i^{\rm T}\frac{\partial V_i^a(x)}{\partial x}, ~i\in \mathcal{V}^a\\
        u_j^p&=-B_{j+m}^{\rm T}\frac{\partial V_j^p(x)}{\partial x}, ~j\in \mathcal{V}^p.
    \end{aligned}
\end{equation}

Next, we show that the Nash equilibrium strategies (\ref{nash}) solve the differential games (\ref{game_a}) and (\ref{game_p}). To prove the asymptotic convergence to $\mathcal{M}_{\rm T}(x^*)$, we first provide the following lemma.

\begin{lemma} \label{lemma4}
    Consider the closed-loop system $\dot x=f(x)+g(y)$ in $\mathbb{R}^{2n}$, if
    \begin{itemize}
        \item[(1)] $g(y)$ tends to 0 asymptotically,
        \item[(2)] $f(x)$ is $\mathcal{R}$-invariant,        
        \item[(3)] there exists a non-trivial equilibrium $x^*$ centered at the origin such that $\frac{\partial f(x^*)}{\partial x}$ has $2n-{\rm dim}(\mathcal{M}_R(x^*))$ eigenvalues with negative real part,
    \end{itemize}
    then the manifold $\mathcal{M}_R(x^*)$ is locally asymptotically stable.
\end{lemma}
\textbf{Proof.} Consider the system $\dot x=f(x)$ in $\mathbb{R}^{2n}$. By Lemma 4.7 of \cite{li2022differential}, $\mathcal{M}_R(x^*)$ is locally asymptotically stable under conditions (2) and (3). Therefore, under conditions (1)-(3), $\mathcal{M}_R(x^*)$ is locally asymptotically stable for the closed-loop system $\dot x=f(x)+g(y)$. $\hfill\square$

\begin{theorem} \label{them2}
    Suppose that $x_i^a(0)\ne x_k^a(0)$ for any $i\ne k,~i,k\in \mathcal{V}^a$, and $x_j^p(0)\ne \frac{1}{m}\sum_{i=1}^m  x_i^a(0)$ for each $j\in \mathcal{V}^p$.
    There exists a $\beta_1^*(x)$ such that for any $\beta_1\le \beta_1^*(x)$, the infinite-horizon differential games in (\ref{game_a}) and (\ref{game_p}) are well defined with $q_i^a(x)$ and $q_j^p(x)$ locally positive semi-definite at $\mathcal{M}_{\rm T}(x^*)$. 
    For the desired manifold $\mathcal{M}_{\rm T}(x^*)$ of two concentric circles, with the designed topology in Section~\ref{sec3.2}, 
    the following value functions
    \begin{align}
        V_i^a(x)= & \frac{\alpha_1}{2n}\sum_{i=1}^m\sum_{j=1}^n \|x_i^a - x_j^p\|^2 \notag\\
        & + \frac{\alpha_2}{d_i^a}\sum_{i=1}^m\sum_{k\in \mathcal{N}_i^a} \frac{1}{\|x_i^a-x_k^a\|}  \label{mm12}\\
        V_j^p(x)= & \beta_1\sum_{j=1}^n  \frac{1}{\|x_j^p-\frac{1}{m}\sum_{i=1}^m x_i^a\|} \notag\\
        & + \frac{\beta_2}{6}\sum_{j=1}^n \sum_{k\in \mathcal{N}_j^p} \|x_j^p-x_k^p\|^3 \label{mm13}
    \end{align}
    for $i\in \mathcal{V}^a$, $j\in \mathcal{V}^p$ are locally positive semi-definite near any submanifold of $\mathcal{M}_{\rm T}(x^*)$, and they also form solutions to the coupled HJB equations (\ref{HJBi}) (\ref{HJBj}) with 
    \begin{equation}
    \begin{aligned}
        Q_i^a&=B_iB_i^{\rm T}+2\sum_{k=1,k\ne i}^m B_kB_k^{\rm T},\\
        Q_j^p&=B_{j+m}B_{j+m}^{\rm T}+2\sum_{k=1,k\ne j}^n B_{k+m}B_{k+m}^{\rm T}.
    \end{aligned}
    \end{equation}

    Furthermore, the feedback controllers in (\ref{nash}) form Nash equilibrium strategies for the games (\ref{game_a}) and (\ref{game_p}), and $x$ converges to $\mathcal{M}_{\rm T}(x^*)$ asymptotically. The radii of the circles are given by Table~\ref{table1}, and the positions of each player and their neighbors follow the result in Lemma~\ref{lemma3}. 
\end{theorem}
\textbf{Proof.} 
Given $C^a(x)$ and $C^p(x)$ in (\ref{mm11}), by (\ref{mm10}) and (\ref{mm16}), direct calculation gives
\begin{equation}
\begin{aligned}
    &q_i^a(x)=\sum_{i=1}^m \left\|\alpha_1(x_c^p-x_i^a)+\alpha_2\sum_{k\in \mathcal{N}_i^a}\frac{x_i^a-x_k^a}{\|x_i^a-x_k^a \|^3} \right \|^2\\
    &\qquad + \frac{2m\alpha_1}{n}\sum_{j=1}^n\!\Big(\!-\beta_1\|x_j^p-x_c^a\|^{-1}\!\\
    &\qquad +\!\beta_2\!\sum_{k\in \mathcal{N}_j^p}\!\|x_j^p\!-\!x_k^p\|(x_j^{p \rm T}-x_c^{a\rm T})(x_j^p\!-x_c^a-\!x_k^p+x_c^a)\Big),\\
    &q_j^p(x)=\sum_{j=1}^n \left \|\beta_1 \frac{x_j^p-x_c^a}{\|x_j^p-x_c^a\|^3} +\beta_2\sum_{k\in \mathcal{N}_j^p}\|x_j^p-x_k^p\|(x_k^p-x_j^p) \right\|^2\\
    &\qquad -\frac{2\alpha_1\beta_1}{n}\sum_{j=1}^n\frac{x_j^{p\rm T}-x_c^{a \rm T}}{\|x_j^p-x_c^a\|}\sum_{j=1}^n(x_j^p-x_c^a).
\end{aligned}
\end{equation}
Since $(x_j^{p \rm T}-x_c^{a\rm T})(x_j^p-x_c^a-x_k^p+x_c^a)\ge 0$, we can always find a $\beta_{11}^*(x)$ such that for any  $\beta_1\le \beta_{11}^*(x)$, it holds that $q_i^a(x)\ge 0$ and $q_j^p(x)\ge 0$ for all $x$. Therefore, the infinite-time differential games in (\ref{game_a}) and (\ref{game_p}) are well defined. 
Let the individual cost $l_i^a(x,u_i^a)=\frac{1}{2}(q_i^a(x)+\|u_i^a\|^2)$ and $l_j^p(x,u_j^p)=\frac{1}{2}(q_j^p(x)+\|u_j^p\|^2)$. It is evident that $l_i^a(x,u_i^a)\ge 0$ and $l_j^p(x,u_j^p)\ge 0$ for any $(x,u)$, and that $l_i^a(x,u_i^a)>0$ and $l_j^p(x,u_j^p)>0$ for any $u_i^a \ne 0$, $u_j^p \ne 0$.

Denote $\mathcal{M}_o=\{x\in \mathbb{R}^{2(m+n)}:C^a(x)x=0,C^p(x)x=0\}$. Then at the desired configuration, we have $\mathcal{M}_{\rm T}(x^*)\subset \mathcal{M}_o$. It is easy to check that the $\mathcal{R}$-invariant and translation invariant functions $V_i^a(x)$ and $V_j^p(x)$ in (\ref{mm12}) and (\ref{mm13}) are solutions to the coupled HJB equations (\ref{HJBi}) and (\ref{HJBj}). By checking an arbitrary coordinate $x^*\in \mathcal{M}_{\rm T}(x^*)$ with radius given in Table~\ref{table1}, based on Lemma~\ref{lemma3}, we obtain $\frac{\partial V_i^a(x^*)}{\partial x}=0$ and $\frac{\partial V_j^p(x^*)}{\partial x}=0$, and each Hessian $\frac{\partial^2 V_i^a(x^*)}{\partial x^2}\ge 0$ and $\frac{\partial^2 V_j^p(x^*)}{\partial x^2}\ge 0$ has three zero eigenvalues on the manifold $\mathcal{M}_{\rm T}(x^*)$. Thus, by Theorem~\ref{them1}, $V_i^a(x)$ and $V_j^p(x)$ both obtain local minimum values at $\mathcal{M}_{\rm T}(x^*)$.
    
Next, we show that $l_i^a(x,0)=0$ and $l_j^p(x,0)=0$ if and only if $x\in \mathcal{M}_{\rm T}(x^*)$. Since $\mathcal{M}_{\rm T}(x^*)\subset \mathcal{M}_o$, we only need to show the necessity. To this end, we note that $l_i^a(x,0)=0$ and $l_j^p(x,0)=0$ if and only if $C^a(x)x=0$ and $C^p(x)x=0$. Define $h^a(x)=C^a(x)x$ and $h^p(x)=C^p(x)x$. It is obvious that $h^a(x)=\frac{\partial V_i^a(x)}{\partial x}$ and $h^p(x)=\frac{\partial V_j^p(x)}{\partial x}$. For any $x\in \mathcal{M}_\mathcal{\rm T}$ and increment $\Delta x$, it follows that 
\begin{equation*}
    h^a(x+b\Delta x)=b\frac{\partial^2 V_i^a(x)}{\partial x^2}\Delta x + o(b),
\end{equation*}
which leads to 
\begin{equation*}
    \|h^b(x+a\Delta x)\|^2=b^2\|\frac{\partial^2 V_i^a(x)}{\partial x^2}\Delta x\|^2+o(b^2).
\end{equation*}
According to Theorem~\ref{them1}, we have
\begin{equation*}
    \Delta x^{\rm T}\frac{\partial^2 V_i^a(x)^{\rm T}}{\partial x^2}\frac{\partial^2 V_i^a(x)}{\partial x^2} \Delta x>0, ~\forall \Delta x\in \mathbb{R}^{2(m+n)}\backslash \mathcal{T}_{\rm T}(x^*).
\end{equation*}
Thus, there must exists a small neighborhood of $\mathbb{R}^{2(m+n)}\backslash \mathcal{T}_{\rm T}(x^*)$ in which $ \|h^a(x+b\Delta x)\|^2>0$ for any $\Delta x\in \mathbb{R}^{2(m+n)}\backslash \mathcal{T}_{\rm T}(x^*)$. Similarly, we can show that $\|h^p(x+b\Delta x)\|^2>0$ for any $\Delta x\in \mathbb{R}^{2(m+n)}\backslash \mathcal{T}_{\rm T}(x^*)$. Since $h^a(x)=0$ and $h^p(x)=0$ in $\mathcal{M}_{\rm T}(x^*)$, the necessity is proved.

The costs $l_i^a(x,u_i^a)$ and $l_j^p(x,u_j^p)$ correspond to the objective function $g(x,u)$ in Lemma~4.5 of \cite{li2022differential} (which is also presented in appendix). It can be seen that all conditions required for $l_i^a(x,u_i^a)$ and $l_j^p(x,u_j^p)$ are satisfied. Thus, by Lemma~4.5 of \cite{li2022differential}, the closed-loop system converges to $\mathcal{M}_{\rm T}(x^*)$ as $t\to \infty$, and the optimal cost-to-go functions are $V_i^a(x_0)=\int_0^{\infty}l_i^a(x,u_i^{a*})dt$ and $V_j^p(x_0)=\int_0^{\infty}l_j^p(x,u_j^{p*})dt$. Since $V_i^a(x_0)$ and $V_j^p(x_0)$ are bounded, we have $x_i^a(t)\ne x_k^a(t)$ for any $i\ne k,~i,k\in \mathcal{V}^a$, and $x_j^p(t)\ne x_c^a(t)$ for each $j\in \mathcal{V}^p$.
Since $u_i^a=-B_i^{\rm T}\frac{\partial V_i^a(x)}{\partial x}=-B_i^{\rm T}C^a(x)x$, $u_j^p=-B_{j+m}^{\rm T}\frac{\partial V_j^p(x)}{\partial x}=-B_{j+m}^{\rm T}C^p(x)x$, and $C^a(x)=0$, $C^p(x)=0$ at the desired manifold, we have $u_i^a=0$ and $u_j^p=0$ at the desired configuration. Therefore, the configuration remains unchanged. 

Next, we claim that $\mathcal{M}_{\rm T}(x^*)$ is asymptotically stable. Denote $x_c^a=\frac{1}{m}\sum_{i=1}^m x_i^a$ and $x_c^p=\frac{1}{m}\sum_{j=1}^n x_j^p$ which are the centers of the active and the passive group, respectively.
Define $\tilde x_i^a=x_i^a-x_c^p$ and $\tilde x_j^p=x_j^p-x_c^a$ for $i=1,\cdots,m$ and $j=1,\cdots,n$, and $\bar x^a=\frac{1}{m}\sum_{i=1}^m\tilde x_i^a$, $\bar x^p=\frac{1}{n}\sum_{j=1}^n\tilde x_j^p$. It is easy to see that $\bar x^a=x_c^a-x_c^p$ and $\bar x^p=-\bar x^a$. 
Under the Nash equilibrium strategies in (\ref{nash}), we have
    \begin{align}
        u_i^a&=\alpha_1(x_c^p-x_i^a)+\alpha_2\sum_{k\in \mathcal{N}_i^a}\frac{x_i^a-x_k^a}{\|x_i^a-x_k^a\|^3} \label{mm14}\\
        u_j^p&=\beta_1 \frac{x_j^p-x_c^a}{\|x_j^p-x_c^a\|^3} +\beta_2\sum_{k\in \mathcal{N}_j^p}\|x_j^p-x_k^p\|(x_k^p-x_j^p).\label{mm15}
    \end{align}
Taking (\ref{mm14}) and (\ref{mm15}) into the dynamics of $\bar x^p$ gives
\begin{equation*}
    \dot {\bar x}^p=-\alpha_1\bar x^p+\beta_1\frac{1}{n}\sum_{j=1}^n\frac{\tilde x_j^p}{\|\tilde x_j^p\|^3}.
\end{equation*}
Let $d=\min\{\|\tilde x_j^p\|\}$, we further have
\begin{equation*}
    \dot {\bar x}^p\le (-\alpha_1+\frac{\beta_1}{d^3})\bar x^p.
\end{equation*}
Choosing $\beta_{12}^*=\alpha_1 d^3$, then $\bar x^p=x_c^p-x_c^a$ converges to 0 exponentially when $\beta_1< \beta_{12}^*$. Then, choosing $\beta_1^*(x)=\min\{\beta_{11}^*,\beta_{12}^*\}$ in Theorem~\ref{them2}, the games are well defined and $\bar x^p=x_c^p-x_c^a$ converges to 0 exponentially.

Since $V_i^a(x)$ is translation invariant, we have $\frac{\partial V_i^a(x)}{\partial x}=\frac{\partial V_i^a(x+x_c^p(\infty))}{\partial x}$, which leads to $u_i^a=-B_i^{\rm T}\frac{\partial V_i^a(x+x_c^p(\infty))}{\partial x}$. Denote $x^a=[x_1^{a \rm T},\cdots,x_m^{a \rm T}]^{\rm T}$ and $x^p=[x_1^{p \rm T},\cdots,x_n^{p \rm T}]^{\rm T}$, we have $\dot x^a=f(x^a)+g(x_c^p-x_c^p(\infty))$. With the controllers $\{u_i^a\}_{i=1}^m$, we can observe that the three conditions of Lemma~\ref{lemma4} are satisfied. Thus, $\mathcal{M}_R(x^{a*})$ is asymptotically stable. Similarly, $\mathcal{M}_R(x^{p*})$ is also asymptotically stable. These facts lead to the conclusion that $\mathcal{M}_{\rm T}$ is asymptotically stable.
$\hfill\square$

\newtheorem{remark}{Remark}
\begin{remark}
    Note that the controller $u_i^a$ in (\ref{mm14}) consists of two parts: the attractive force from the passive group that steers the active players move towards the passive players, and the repelling force among active members to maintain separation. Similarly, the controller $u_j^p$ consists of repelling forces, driving the passive players to flee, and attractive forces from neighbors, ensuring the formation of the desired configuration. 
\end{remark}

The controllers drive the active players to converge to a circle with radius $r_a$, while the passive players converge to other circle with radius $r_p$. If $r_a>r_p$, the active group achieves the encirclement control. Conversely, if $r_p>r_a$, the passive group attains counter-encirclement control. Besides, the containment control can also be achieved as a result of the games. Before presenting the result, we give the definition of convex hull.
\begin{definition}
    A set $\mathcal{C}\subseteq \mathbb{R}^n$ is convex if $(1-\rho)x+\rho y\in \mathcal{C}$ for any $x,y\in \mathcal{C}$ and any $\rho \in [0,1]$. The convex hull ${\rm Co}(X)$ of a finite set of point $X=\{x_1,x_2,\cdots,x_n\}$ is the minimal convex set containing all points in $X$.
\end{definition} 
\begin{lemma}
        (Containment control) Suppose that the players asymptotically converge to the desired manifold under the controllers (\ref{mm14}) and (\ref{mm15}). When $m\ge 3,~n\ge 2$ or $m,n=2$, if $r_p\le r_acos\frac{\pi}{m}$, then the passive players asymptotically converge to the convex hull formed by the active group ; when $m\ge 2,~n\ge 3$ or $m,n=2$, if $r_a\le r_pcos\frac{\pi}{n}$, then the active players asymptotically converge to the convex hull formed by the passive group.
\end{lemma}
\textbf{Proof.} According to the definition, on the desired circle, the convex hull formed by the active group is the set containing all points within and on the $m$-sided polygon formed by the active players. The shortest distance from the center of the circle to any points on the polygon is given by $r_acos\frac{\pi}{m}$. Therefore, if $r_p\le r_acos\frac{\pi}{m}$, the passive players lie in the convex hull formed by the active group. Similarly, the active players are in the convex hull formed by the passive group if $r_a\le r_pcos\frac{\pi}{n}$.


\section{Numerical simulation}\label{sec5}
In this section, we present numerical simulations to demonstrate our encirclement control results under the designed undirected topologies in Section~\ref{sec3.2} and Nash equilibrium strategies (\ref{nash}). 

\begin{figure}	
	\centering	\includegraphics[width=0.45\textwidth]{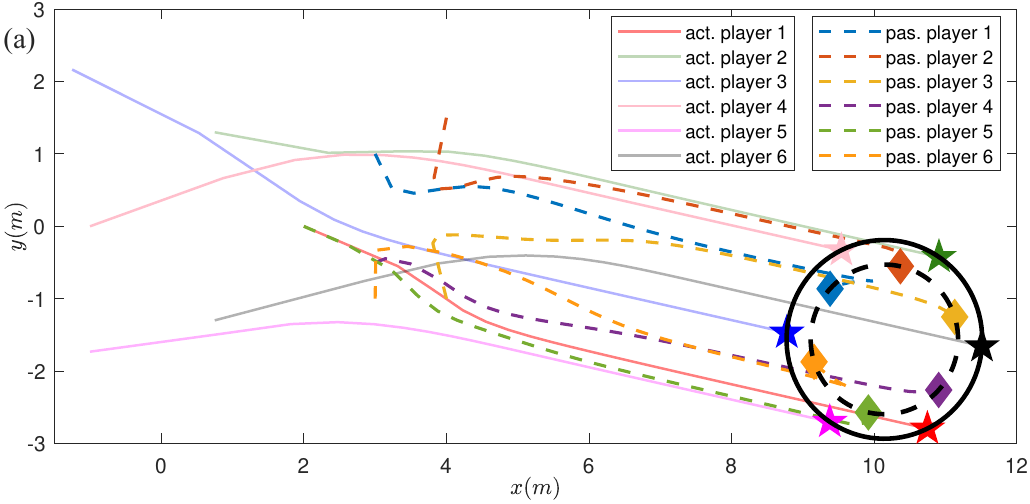}
    \includegraphics[width=0.45\textwidth]{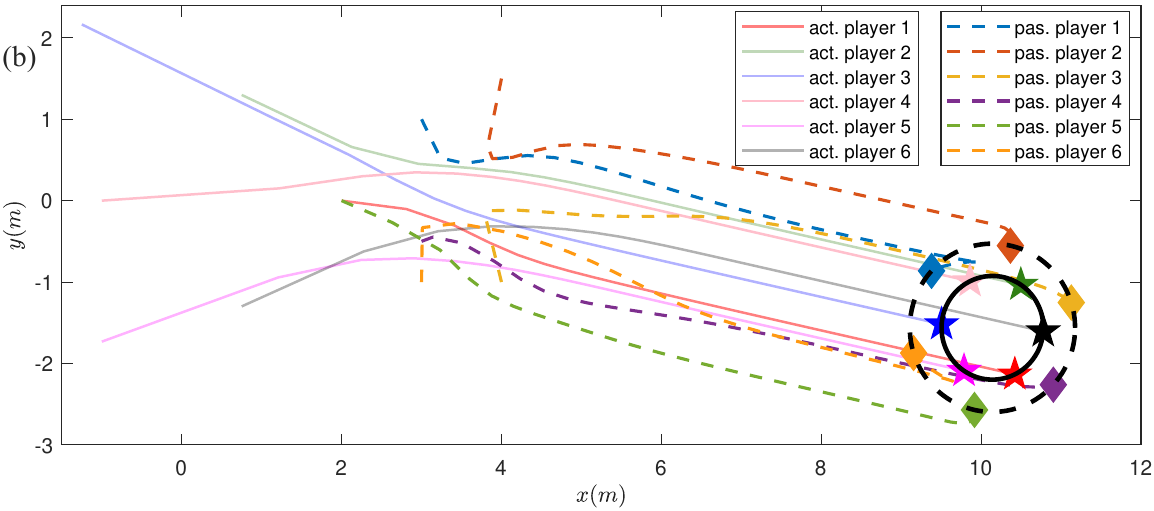}
    \includegraphics[width=0.19\textwidth]{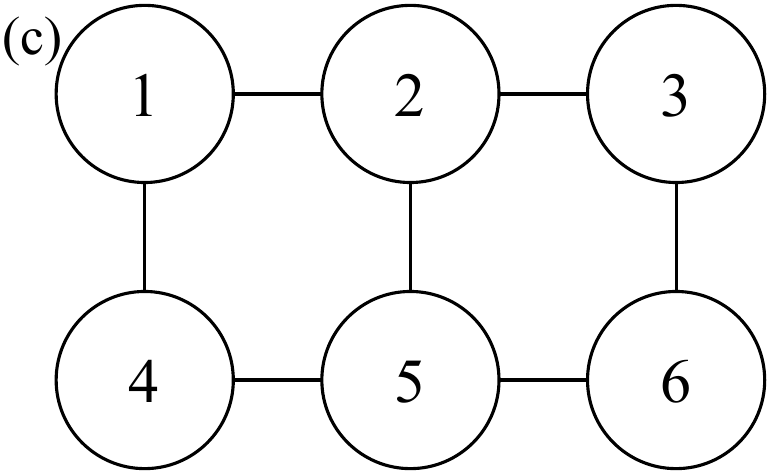}~~~~
     \includegraphics[width=0.19\textwidth]{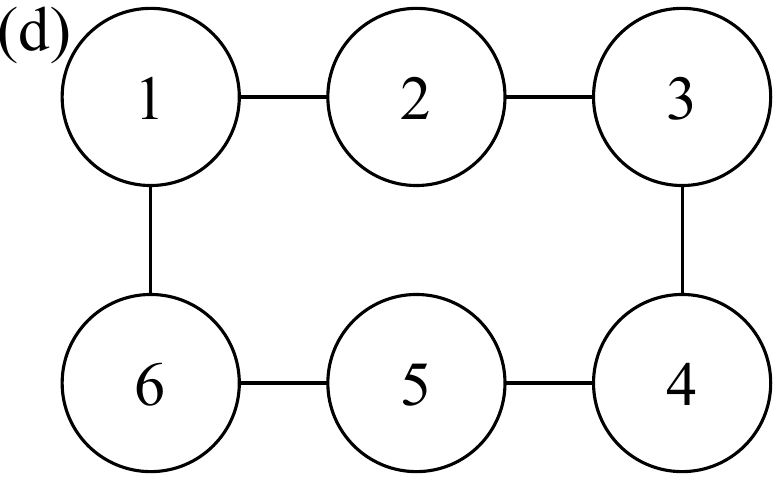}
	\caption{Results of the case $m=6$, $n=6$. (a) The active group achieves encirclement control under the settings $\alpha_1=16$, $\alpha_2=50$, $\beta_1=8$, $\beta_2=7$. (b) The passive group achieves the counter-encirclement control under the settings $\alpha_1=16$, $\alpha_2=5$, $\beta_1=8$, $\beta_2=7$. (c) Topology of the active group. (d) Topology of the passive group.}
	\label{case_1}
\end{figure}

\begin{figure}	
	\centering	\includegraphics[width=0.45\textwidth]{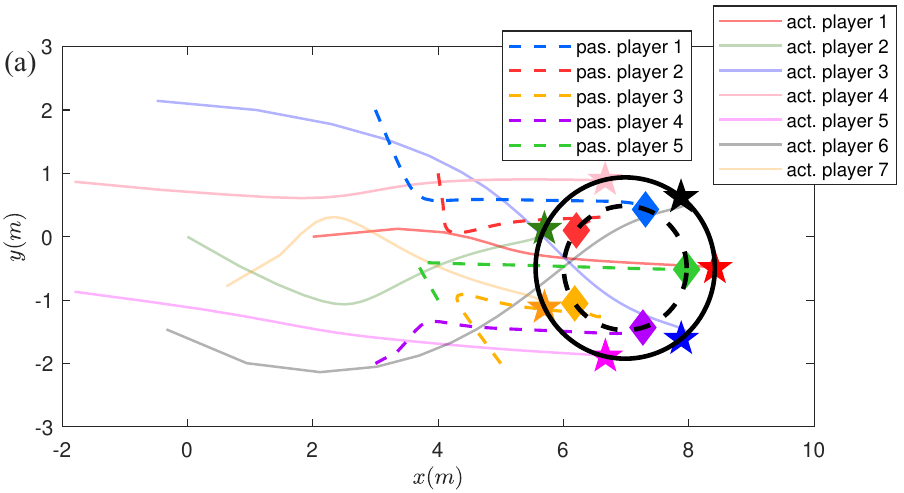}
    \includegraphics[width=0.45\textwidth]{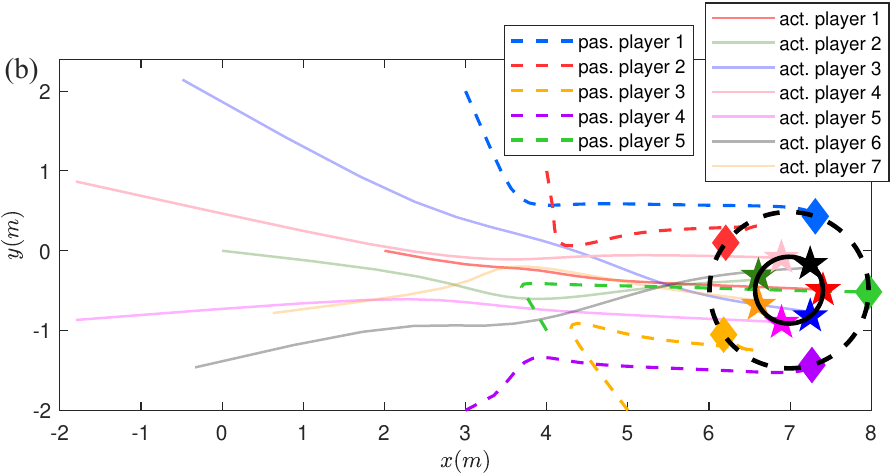}
    \includegraphics[width=0.24\textwidth]{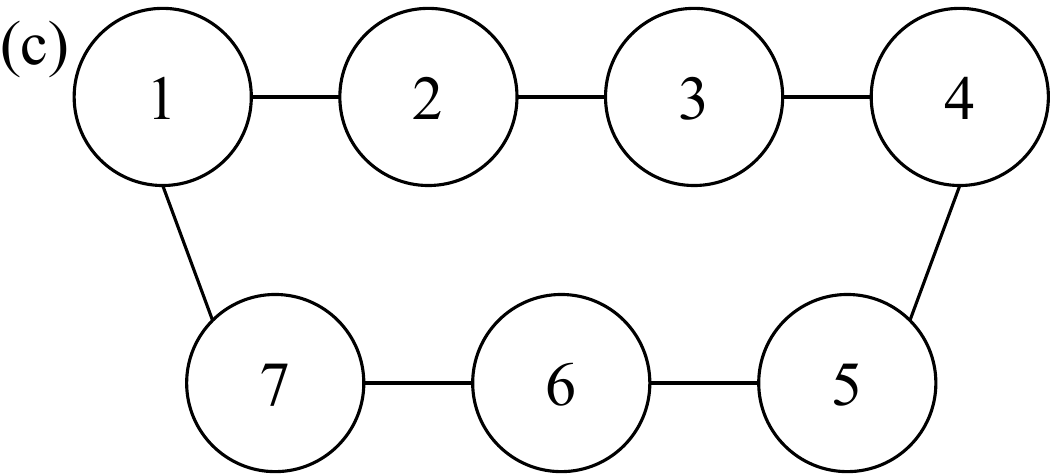}~~~
     \includegraphics[width=0.18\textwidth]{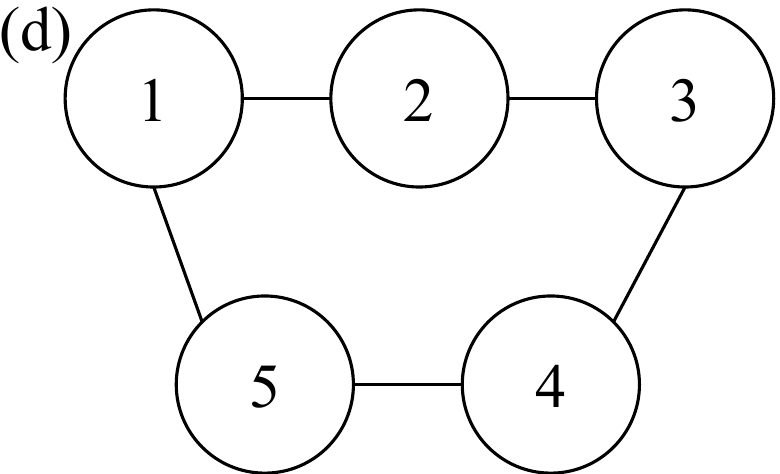}
	\caption{Results of the case $m=7$, $n=5$. (a) The active group achieves encirclement control under the settings $\alpha_1=7$, $\alpha_2=40$, $\beta_1=3$, $\beta_2=2$. (b) The passive group achieves the counter-encirclement control under the settings $\alpha_1=7$, $\alpha_2=1$, $\beta_1=3$, $\beta_2=2$. (c) Topology of the active group. (d) Topology of the passive group.}
	\label{case_2}
\end{figure}

\begin{figure}	
	\centering	\includegraphics[width=0.47\textwidth]{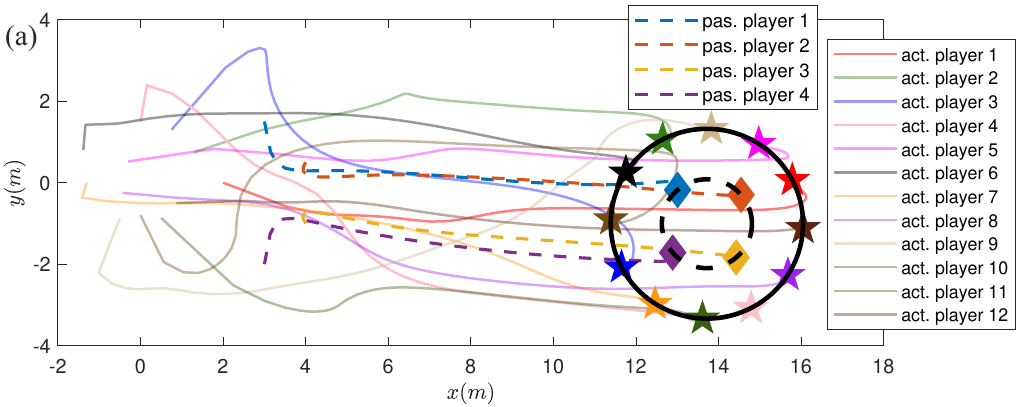}
    \includegraphics[width=0.47\textwidth]{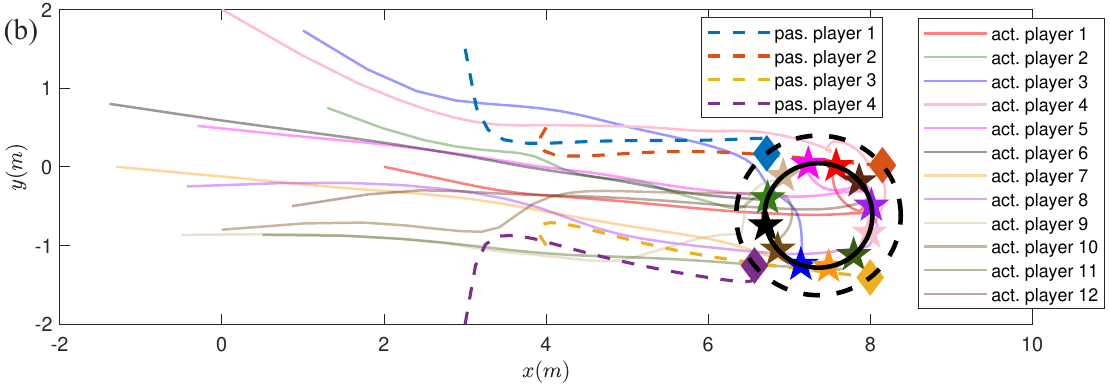}
    \includegraphics[width=0.3\textwidth]{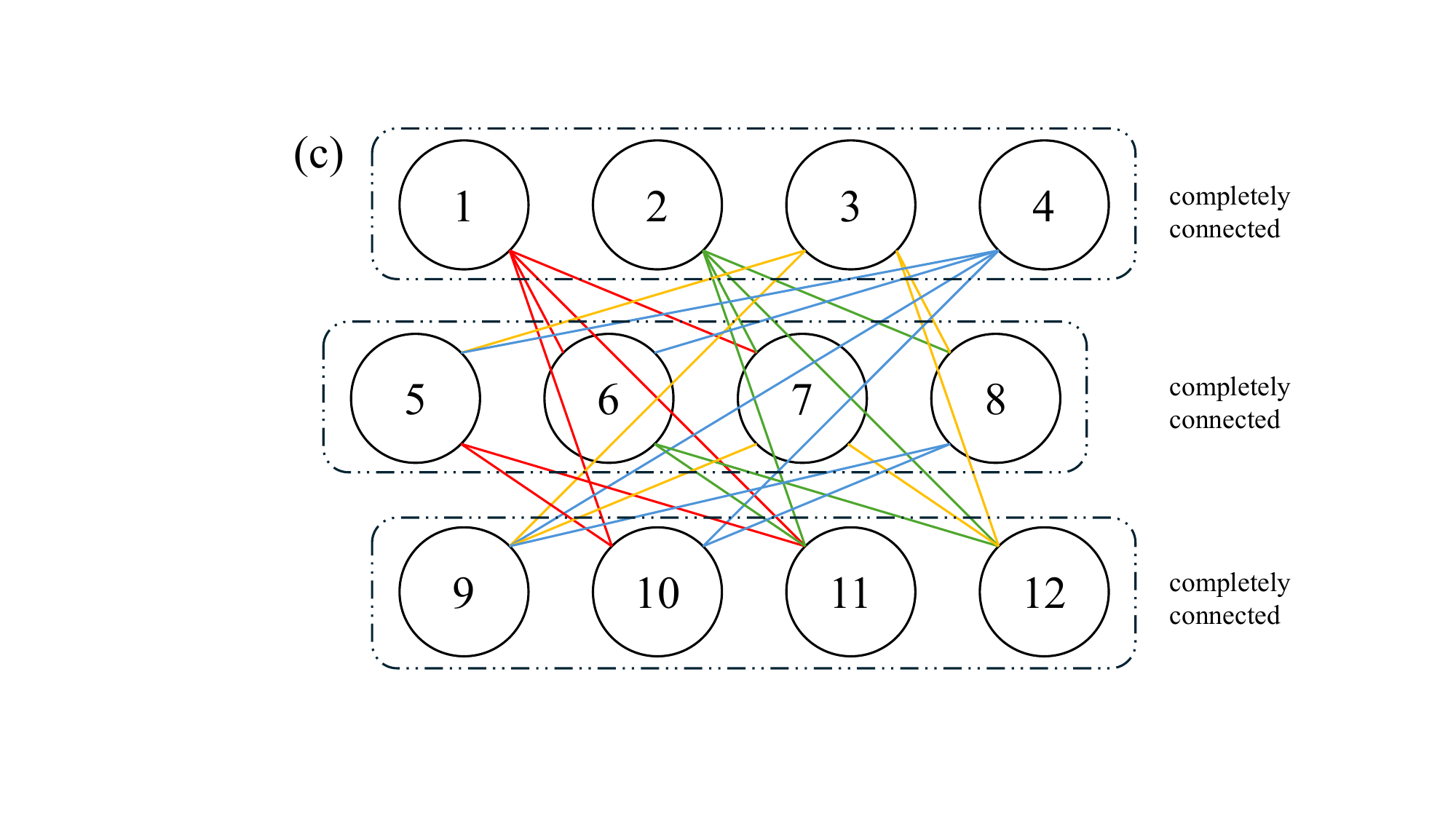}~~~~
     \includegraphics[width=0.13\textwidth]{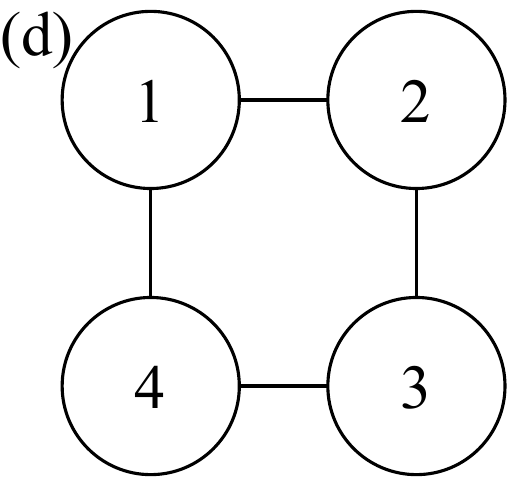}
	\caption{Results of the case $m=12$, $n=4$. (a) The active group achieves encirclement control under the settings $\alpha_1=6.5$, $\alpha_2=40$, $\beta_1=8$, $\beta_2=2$. (b) The passive group achieves the counter-encirclement control under the settings $\alpha_1=7$, $\alpha_2=1$, $\beta_1=6$, $\beta_2=2$. (c) Topology of the active group. (d) Topology of the passive group.}
	\label{case_3}
\end{figure}

The simulation results are shown in Figs.~\ref{case_1}--\ref{case_3}, where the initial positions of players are randomly given, and the topology of the passive group is an undirected ring. These examples cover the three cases in Assumption~\ref{ass1} for the active group, which correspond to Fig.~\ref{case_1}, Fig.~\ref{case_2}, and Fig.~\ref{case_3}, respectively. The results demonstrate that the Nash equilibrium strategies successfully drive the players toward the desired relative configuration, and both encirclement and counter-encirclement behaviors are achieved under different parameters settings for $\alpha_1$, $\alpha_2$, $\beta_1$, and $\beta_2$.   
It is also worth noting that for each active player, its neighbors are located at the farthest position on its circle, while for each passive player, its neighbors are at the nearest point on its respective circle.

\section{Conclusion}\label{sec6}
In this paper, we have investigated the intrinsic encirclement control problem by formulating it as an infinite-time differential game. Different from most literature that predefine an expected radius or desired path around the target(s), our approach drives agents to converge to the desired formation manifold, by allowing permutation, rotation, and translation of players. The Nash equilibrium strategies are obtained in an intrinsic way in the sense that they depends only on the inter-agent interactions and geometric properties of the network. Our results demonstrate that both encirclement control by the active group and counter-encirclement control by the passive group can be achieved. In the near future, we plan to extend our results to higher-order systems with directed interaction graphs.

\section*{Appendix}

\begin{figure*}
\begin{equation} \label{mm5}
\begin{aligned}
    Hv_{k,l} &=\begin{bmatrix}
        (C_0 + w^{k(m-1)}RC_{m-1} + w^{k(m-2)}R^2C_{m-2} + \cdots + w^kR^{m-1}C_1)\psi_j \\
        (C_1 + w^{k(m-1)}RC_0 + w^{k(m-2)}R^2C_{m-1} + \cdots + w^kR^{m-1}C_2)\psi_j\\
        \vdots\\
        (C_{m-1} + w^{k(m-1)}RC_{m-2} + w^{k(m-2)}R^2C_{m-3} + \cdots + w^kR^{m-1}C_0)\psi_j
    \end{bmatrix}\\
    &=\begin{bmatrix}
        (C_0 + w^{k(m-1)}RC_{m-1} + w^{k(m-2)}R^2C_{m-2} + \cdots + w^kR^{m-1}C_1)\psi_j \\
        (w^kC_1 + RC_0 + w^{k(m-1)}R^2C_{m-1} +\cdots + w^{2k}R^{m-1}C_2)w^{k(m-1)}R^{-1}R\psi_j\\
        \vdots\\
        (w^{k(m-1)}C_{m-1} + w^{k(m-2)}RC_{m-2} + w^{k(m-3)}R^2C_{m-3} + \cdots + R^{m-1}C_0)w^kR^{1-m}R^{m-1}\psi_j
    \end{bmatrix}\\
    &=\begin{bmatrix}
        (C_0 + w^{k(m-1)}RC_{m-1} + w^{k(m-2)}R^2C_{m-2} + \cdots + w^kR^{m-1}C_1)\psi_j \\
        R^{m-1}(w^kC_1 + RC_0 + w^{k(m-1)}R^2C_{m-1} + \cdots + w^{2k}R^{m-1}C_2)R^{1-m}w^{k(m-1)}R^{-1}R\psi_j\\
        \vdots\\
        R(w^{k(m-1)}C_{m-1} + w^{k(m-2)}RC_{m-2} + w^{k(m-3)}R^2C_{m-3} + \cdots + R^{m-1}C_0)R^{-1
        }w^kR^{1-m}R^{m-1}\psi_j
    \end{bmatrix}\\
      &=\begin{bmatrix}
        (C_0 + w^{k(m-1)}RC_{m-1} + w^{k(m-2)}R^2C_{m-2} + \cdots + w^kR^{m-1}C_1)\psi_j \\
        (C_0 + w^{k(m-1)}RC_{m-1} + w^{k(m-2)}R^2C_{m-2} + \cdots + w^kR^{m-1}C_1)w^{k(m-1)}R\psi_j\\
        \vdots\\
        (C_0 + w^{k(m-1)}RC_{m-1} + w^{k(m-2)}R^2C_{m-2} + \cdots + w^kR^{m-1}C_1)w^kR^{m-1}\psi_j
    \end{bmatrix}\\
    & = \lambda_l(D_k)v_{k,l},~~l=1,2.
\end{aligned}
\end{equation}
\end{figure*}

\textbf{Proof of Lemma~\ref{lemma2}.} As $H$ is symmetric and each block $C_i$, $i=0,1,\cdots,m-1$ is also symmetric, $H$ can be rewritten as 
\begin{equation*}
    H\!\!=\!\!\left[\begin{array}{ccccc}
C_{0} & R C_{m\!-\!1} R^{-1} & R^{2} C_{m-2} R^{-2} & \cdots & R^{m-1} C_{1} R^{1-m} \\
C_{1} & R C_{0} R^{-1} & R^{2} C_{m-1} R^{-2} & \cdots & R^{m-1} C_{2} R^{1-m} \\
C_{2} & R C_{1} R^{-1} & R^{2} C_{0} R^{-2} & \cdots & R^{m-1} C_{3} R^{1-m} \\
\vdots & \vdots & \vdots & \cdots & \vdots \\
C_{m-1} & R C_{m-2} R^{-1} & R^{2} C_{m-3} R^{-2} & \cdots & R^{m-1} C_{0} R^{1-m}
\end{array}\right]\!.
\end{equation*}


By direct calculating, we have equation (\ref{mm5}). It implies that the eigenvalues of $H$ consist of the eigenvalues of $B_k$, and the corresponding eigenvectors are $v_{k,l}$. $\hfill\square$

\textbf{Lemma~4.5 of \cite{li2022differential}.}
    Consider an infinite time optimal control problem
    \begin{equation} \label{opm}
    \begin{aligned}
        &\min ~\int_0^{\infty} g(x,u,\theta(x))dt \\
        & s.t.~\dot x=f(x,u),\\
        &\quad ~~~ x(0)=x_0,
    \end{aligned} 
    \end{equation}
    where $g$ and $f$ are Lipschitz continuous functions, $\theta(x)=[\theta_1(x),\cdots,\theta_l(x)]^{\rm T}$ is a set of parameters. 

    Denote $\mathcal{M}$ as a closed manifold. Assume the following assumptions hold:
    \begin{itemize}
        \item $g(x,0,\theta(x))=0$ if $x\in \mathcal{M}$.
        \item For some $\theta(x)$, $g(x,u,\theta(x))\ge 0$ for any $(x,u)$, and $g(x,u,\theta(x))>0$ if $u\ne0$.
        \item  $g(x(t),0,\theta(x))=0$ for all $t\ge 0$ if and only if $x(t)\in \mathcal{M}$ for all $t\ge 0$.
    \end{itemize}
    If there exists a real-valued function $V\in \mathbb{C}^1$ that is positive definite with respect to $\mathcal{M}$, that is, $V(x)\ge 0$ for any $x$, and $V(x)=0$ if $x\in \mathcal{M}$, and satisfies
    \begin{equation*}
        \min_{u}\{g(x,u,\theta(x))+\frac{\partial V(x)}{\partial x}f(x,u)\}=0,
    \end{equation*}
    then
    \begin{itemize}
        \item  the controller        \begin{equation*}
            u^*={\arg \min}_u\{g(x,u,\theta(x))+\frac{\partial V(x)}{\partial x}f(x,u)\}
        \end{equation*}
        is the optimal one to (\ref{opm}) and the closed-loop system converges to $\mathcal{M}$ as $t\to \infty$.
        \item $V(x_0)=\min_{u^*} \int_0^{\infty}g(x,u,\theta(x))dt$ is the optimal cost-to-go function.
    \end{itemize}


\end{document}